\documentclass[11pt,a4paper]{article}
\usepackage{amsmath,amssymb}
\usepackage{amsmath}

\usepackage{color}
\usepackage[colorlinks,linkcolor=blue,citecolor=red]{hyperref}

\DeclareMathSymbol{\bbbr}{\mathalpha}{AMSb}{"52}
\DeclareMathSymbol{\bbbc}{\mathalpha}{AMSb}{"52}

\newtheorem{theorem}{Theorem}

\newtheorem{corollary}[theorem]{Corollary}

\newtheorem{definition}[theorem]{Definition}

\newtheorem{lemma}[theorem]{Lemma}

\newtheorem{proposition}[theorem]{Proposition}

\begin{document}

\title{ Darboux integrability for diagonal systems of hydrodynamic type}

\author{{\Large Sergey I. Agafonov}\\
\\
Department of Mathematics,\\
S\~ao Paulo State University-UNESP,\\ S\~ao Jos\'e do Rio Preto, Brazil\\
e-mail: {\tt sergey.agafonov@gmail.com} }
\date{}
\maketitle
\unitlength=1mm

\vspace{1cm}

\begin{abstract} We prove that 1) diagonal systems of hydrodynamic type are Darboux integrable if and only if the corresponding systems for commuting flows are Darboux integrable,  2) systems for commuting flows are Darboux integrable if and only if  the Laplace transformation sequences  terminate, 3) Darboux integrable systems are necessarily semihamiltonian. We give geometric interpretation for Darboux integrability of such systems in terms of congruences of lines and in terms of solution orbits with respect to symmetry subalgebras, discuss known and new examples.
\bigskip


\bigskip

\end{abstract}

%\vspace{-7mm}
%\newpage

%\tableofcontents

\section{Introduction} The notion of integrability for partial differential equations evolved dramatically over time: from "explicit" solutions,  like the d'Alambert formula for one-dimensional  wave equation, to highly sophisticated techniques such as inverse scattaring transform for KdV and Schr\"odinger equations.  One of the classical approaches to integrability, worked out by Monge with "intermediate integrals" and further  developed by Darboux, to deal with scalar hyperbolic equations of 2nd order with two independent variables, is adaptable to a larger class of systems.
In this paper,  we apply this approach to diagonal systems of hydrodynamic type
\begin{equation}\label{hydro}
u^i_t=\lambda^i(u)u^i_x, \ \ \ i=1,...,n,\ \ \ \mbox{(no summation).}
\end{equation}
Such systems are important  not only for studying classical  \cite{R-60} and modern integrable models \cite{DN-83,T-00} of mathematical physics, they also turned out  to be a key tool for exploring geometric structures (see, for example, \cite{Z-98,BM-11,A-21}).

The Darboux method has its roots in the geometric theory of partial differential equations \cite{S-00,IL-03}, whose foundations were laid by S.Lie, E.Vessiot, and \'E. Cartan. It works for hyperbolic equations such that each its {\it Monge characteristic distribution} has at least two independent {\it Riemann invariants}.

System (\ref{hydro}) is automatically hyperbolic and each its Monge  distribution  has at least one Riemann invariant, namely $u^i$.

\begin{definition}
System (\ref{hydro})  is Darboux integrable if  each Monge characteristic distribution of the system or of  its prolongation admits at least one extra Riemann invariant.
\end{definition}
If the system is Darboux integrable then solving it reduces to integration of a complete exterior differential system, i.e. to integration of ordinary differential equations.

It turns out that the Darboux integrability criterion can be written in terms of the functions
\begin{equation}\label{aij}
a_{ij}=\frac{\partial_j\lambda^ i}{\lambda^j-\lambda^i}
\end{equation}
 and their derivatives (we use the notation $\partial _i=\partial_{u^i}$). In particular, Darboux integrable systems are necessary semihamiltonian, i.e.
\begin{equation}\label{semiH}
\partial_ja_{ik}=\partial_ka_{ij},  \  \  \  \ \forall \ i\ne  k\ne j\ne i.
\end{equation}
The first main result of the paper is that this is true if some prolongation is Darboux integrable, no matter how high the order is. 

Semihamiltonian diagonal systems  (\ref{hydro}) have very reach geometric and analytic structure \cite{T-00}. In particular, they have infinitely many conservation laws and, the most remarkably, they  can be linearized by the so-called generalized hodograph method: any solution with $u^i_x\ne 0$ $\forall i$ has a local representation in the following implicit form (see \cite{T-90}):
\begin{equation}\label{Tsarevform}
\mu^i(u)=\lambda^i(u)t+x,\ \ \ \ i=1,...,n,
\end{equation}
where $\mu$ is a suitable solution of the linear system
\begin{equation}\label{commuting}
\partial_k\mu^i=a_{ik}(\mu^k-\mu^i),\ \ \ \ i\ne k, \ \ \ \mbox{(no summation).}
\end{equation}
Conversely, a generic solution $\mu$ to (\ref{commuting}) gives locally a solution  to the semihamiltonian system (\ref{hydro}).
Observe that $\lambda$-s also  solve this system. System (\ref{commuting}) defines an infinite dimensional  commutative algebra of flows, commuting with (\ref{hydro})  \cite{T-90}: equations (\ref{hydro}) and the following ones
\begin{equation}\label{hydroMu}
u^i_{\tau}=\mu^i(u)u^i_x
\end{equation}
are compatible, their common solution is a vector function $u(t,\tau,x)$.

We also study Darboux integrability for the linear system (\ref{commuting}) and show that the $r$-th prolongation of (\ref{hydro}) is Darboux integrable if and only if the $r$-th prolongation of (\ref{commuting}) is Darboux integrable. In particular, this implies that if system (\ref{hydro}) is Darboux integrable then any commuting system (\ref{hydroMu}) is also Darboux integrable.

If for some $i$ the coefficients $a_{ik}$ do not vanish, one can exclude all $\mu^k$, $k\ne i$ from (\ref{commuting}) and obtain a compatible linear system for $\mu^i$. For such systems, Darboux  generalized the Laplace method originally developed for just one linear hyperbolic equation with two independent variables (see \cite{D-96} IV, p.274). We show that Laplace transformations for this system are liftable to all $\mu$ (and therefore to all $\lambda$) and that  the $i$-th Monge distribution has an extra Riemann invariants of order $r$ if and only if  such $i$-sequence of Laplace transformations terminates after $r$ steps (this happens, for example, when  for the transformed system holds $\bar{a}_{im}=0$.). This is completely analogous to the results of \cite{SZ-95,JA-97,AK-98,ZS-01}: a nonlinear 2nd order hyperbolic equation  with two independent variables is Darboux integrable if and only if the Laplace transformation sequences  of its linearizations terminate. (For a linear equations, integrability for systems with such terminating sequences was understood already by Laplace.  Darboux (see \cite{D-96} II, p.23) describes this method though he did not mention any "intermediate integrals"!?) Remarkably, the Laplace transformation for semihamiltonian systems (\ref{hydro}) defined in this paper via linear system for commuting flows (\ref{commuting}) coincides with the one defined via the linear system for the densities of conservation laws of (\ref{hydro}) in \cite{F-97} (see also \cite{KT-97} and \cite{AF-99}). We also generalize Laplace transformaion for the cases when only some of  the coefficients $a_{ik}$ do not vanish. 

The set of Darboux integrable systems (\ref{hydro}) is sufficiently large.  In section \ref{sDi}, we show that, for Darboux integrability by integrals of the first order,  it can be parametrized by $n(n+1)$ functions of one variable and discuss how integrate them by solving ODEs either via generalized hodograph transform or directly constructing complete involutions for (\ref{hydro}). For linearly degenerate systems (i.e. such that $\lambda^i_i=0$), the formulas for such involutions can be written explicitly and   were found in \cite{T-85,P-87}. In fact, they are explicit because the second Riemann invariant can be computed explicitly. We describe a pair of more general classes when this is the case.  We also show that for the Darboux integrable systems (\ref{hydro}) the corresponding system (\ref{commuting}) can be integrated in quadratures.

Semihamiltonian systems can be defined as strictly hyperbolic systems of conservation laws admitting $n$ Riemann invariants \cite{S-94}. Choosing $n$ functionally independent densities $N^i(u)$, we rewrite (\ref{hydro}) as a system of conservation laws
\begin{equation}
\frac{d}{dt}N^i(u)=\frac{d}{dx}M^i(u),~~~i=1,...,n.
\label{CL}
\end{equation}
There is a natural projective geometry associated to a system of conservation laws (\ref{CL}). Namely,  a
 {\it congruence of lines} (see \cite{AF-96}), i.e. an $n$-parameter family of lines
\begin{equation}
Y^i=N^i(u)\ Y^0+M^i(u)Y^{n+1}, ~~~ i=1,..., n,
\label{congrn}
\end{equation}
in $(n+1)$-dimensional projective space $\mathbb P^{n+1}$ with homogeneous coordinates
$[Y^0:...:Y^{n+1}]$ carries all the information about the system.
The basic concepts of the  theory of conservation law systems
(rarefaction curves, shock curves, reciprocal transformations, linear degeneracy, etc.) translate  in geometrical language \cite{AF-96,AF-99}.  In fact, this language extends also to a geometric interpretation of solutions \cite{A-20}.

In section \ref{sectionCongr}, we present geometric counterparts of Laplace transformations and Darboux integrability. For better visualization, here in introduction we exemplify them for $n=2$. Any line of the congruence (\ref{congrn}) touches two so-called {\it focal} surfaces $F_1$ and $F_2$ in $\mathbb{P}^3$. Coordinate lines $u^2=const$ correspond to developable surface of the congruence with the cuspidal edges on $F_1$. That means that the congruence lines are tangent to $F_1$ along curves of some family $C_1$ as only the first coordinate $u^1$ is varying.  Similarly, the congruence lines are  tangent to $F_2$ along curves of some family $C_2$ as only the second coordinate  $u^2$ varies. If $u^1=const$ then  the points, where the lines touches $F_1$,  trace  curves  of a family $K_1$, conjugate to ones in $C_1$. The family $K_2$ on $F_2$ is defined similarly. Laplace transformations of the congruence are defined as changing the role of the families $C_i$ and $K_i$ on one of the focal surfaces: the lines of the transformed congruence are tangent to the curves $K_i$ \cite{F-97,AF-99}. For example, if $i=1$ then one of the focal surfaces  of the new congruence coincides with $F_1$ and the second is formed by the cuspidal edges of the developable surfaces of the new congruence, these surfaces being cut from the congruence by $u^2=const$.  This construction for line congruences in $\mathbb{P}^3$ was given by Darboux (see \cite{D-96}, V2, book 4, cap. I) as a geometric interpretation for Laplace transformation of the equation for the parametrization of the  focal surface by the conjugate coordinates $u^1,u^2$.  We show that this construction provides also the geometric picture of Laplace transformations for the system (\ref{commuting}). Note that the equation for the corresponding $\lambda$ is not the equation for the above parametrization of the focal surface, therefore terminating of a sequence of Laplace transformations for (\ref{commuting}) is geometrically different from the one given by Darboux: the curves, conjugated to the cuspidal edges on the new focal surface, are cut by a pencil of planes, thus making  impossible the next Laplace transformation "starting from" this new surface, while in the Darboux picture  the new focal surface of the transformed congruence degenerates into a curve.

Another geometric aspect, discussed in this paper, is related to the symmetries of (\ref{hydro}).  Darboux integrability, manifesting a kind of variable separation (see \cite{AFV-09}), is reflected also in the structure of the symmetry algebra of commuting flows. Fixing a relation between two Riemann invariants of the same Monge distribution, we get a differential relation compatible with the system and depending on one arbitrary function of one variable. We show that this relation fixes an orbit of a solution with respect to some  subalgebra of the symmetry algebra, defined by (\ref{commuting}).  Moreover, this subalgebra is also cut  by a differential relation, uniquely defined by the choice of Monge distribution.

Even though Darboux method is potentially applicable to systems of equations for vector-valued functions in any number of variables, most of the studies were devoted to a single 2nd order equation for scalar funcion in two variables (see \cite{P-60,R-99,H-01} for generalization of the method to higher order equations  and \cite{T-05} for extension of Laplace transformations on systems of linear hyperbolic PDEs with several dependent and two independent variables via  {\it generalized factorization}).  
 This can be explained by  the increasing complexity of calculations needed to work out higher-dimensional cases. This complexity, probably, explains also the scarcity of explicitly elaborated  examples within prescribed classes. We hope that this paper contributes to a better understanding of the case of hydrodynamic systems.

 A warning on terminology: systems (\ref{hydro}) are "hyperbolic" in the sense of the theory developed in \cite{BGH-95} only for $n=2$.

\section{Monge characteristic distribution for diagonal systems}
Let us rewrite equations (\ref{hydro}) as an exterior differential system in $\mathbb{R}^{2n+2}$ with coordinates $u^i,u^i_x,x,t$:
\begin{equation}\label{ext}
\sigma_i:=u^i_x(dx+\lambda^i(u)dt)-du^i=0.
\end{equation}
Then solutions to (\ref{hydro}) are integral  manifolds of this system with $dx\wedge dt\ne 0$. In the dual language, integration of  (\ref{ext}) amounts to finding 2-dimensional integral manifolds tangent to the distribution $\cal{V}$ with generators
\begin{equation}\label{V}
\mathcal{V}:=\left\langle  \partial_x+\sum_k u^k_x\partial_k,\  \partial_t+\sum_k\lambda^k u^k_x\partial_k,\  \partial_{u^i_x} \right\rangle_{i=1,...,n}
\end{equation}
\begin{definition}\cite{S-00}\label{M}
A subdistribution $\mathcal{M}$ of $\mathcal{V}$ is Monge characteristic distribution if
\begin{enumerate}
\item $\mathcal{M}$ is singular: for any vector $m\in \mathcal{M}$ holds
$$
\dim \langle m\lrcorner d\sigma_i,\ i=1,...,n \rangle|_{\mathcal{V}}<
\max_{v\in \mathcal{V}} \dim \langle v\lrcorner d\sigma_i,\ i=1,...,n \rangle|_{\mathcal{V}},
$$
\item $\mathcal{M}\cap \mathcal{I}\ne 0$ for any maximal involution $\mathcal{I}$ of $\mathcal{V}$,
\item $\mathcal{M}\cap \mathcal{W}$ is complete (i.e. integrable) for any maximal complete subdistribution $\mathcal{W}$ of $\mathcal{M}$.
\end{enumerate}
\end{definition}
(Recall that involution is  a subdistribution $\mathcal{I}$ of $\mathcal{V}$ such that $[\mathcal{I},\mathcal{I}]\in \mathcal{V}$.)

An easy exercise in linear algebra shows that (\ref{hydro}) has $n$ Monge characteristic distributions and the $i$-th one  is cut by the following equations:
\begin{equation}\label{Monge1}
\begin{array}{l}
dx+\lambda^idt=0,\\
du^i=0,\\
du^j+(\lambda^i-\lambda^j)u^j_xdt=0,\  j\ne i \\
du^i_x-D_x(\lambda^i)u^i_xdt=0,
\end{array}
\end{equation}
 where $D_x=\partial_x+\sum_ku^k_x\partial_k+...$ is the total derivative by $x$.
 It is immediate that this distribution has at least one Riemann invariant (i.e. first integral of the Monge characteristic distributions), namely $u^i$.
The system is Darboux integrable if each Monge distribution admits at least one extra Riemann invariant. Darboux inegrability depends on the codimension of  the exterior differential ideal generated by (\ref{Monge1}).

One can consider $r$-th prolongations of system (\ref{hydro}) and rewrite these prolongations also as an exterior differential systems in $\mathbb{R}^{(2+r)n+2}$ with coordinates $u^i,u^i_x,u^i_{xx},..., x,t$. For these "prolonged" exterior systems, Monge characteristic distributions are defined in a similar way.

We say that a Riemann invariant (of a prolonged system) is of order  $q$ if it depends on derivatives of order $\le q$ and is not of order less then $q$.
\begin{lemma}\label{order01} Let the system (\ref{hydro}) be strictly hyperbolic (i.e. $\lambda^i\ne \lambda^j$ for $i\ne j$) then

a) the $i$-th Monge distribution  admits one extra Riemann invariant of order 0  if and only if $$\lambda^i_k=0\ \ \forall \ k\ne i ,$$

b) the $i$-th Monge distribution  admits one Riemann invariant of order 1 and only one Riemann invariant of order 0  if and only if the following conditions hold true:
\begin{itemize}
\item $\partial_j(a_{ik})=\partial_k(a_{ij})$ $\forall \ k\ne j\ne i \ne k,$
\item $\exists \ j\ne i: \  \lambda^i_j\ne 0,\ \ \partial_jb_{ij}+a_{ij}b_{ij}=0$, where $b_{ij}=\frac{\lambda^i_{ij}}{\lambda^i_j}+\frac{\lambda^i_i}{\lambda^j-\lambda^i} = \frac{\partial_i a_{ij}}{a_{ij}}-a_{ji},$
\item if $\lambda ^i_k\ne 0$ for $k\ne i$  then $b_{ik}=b_{ij}$ and  $\partial_kb_{ik}+a_{ik}b_{ik}=0.$
\end{itemize}

\end{lemma}

 {\it Proof:}
The $i$-th Monge distribution is generated by the  vector fields $\partial _{u^k_x}$ with $k\ne i$ and by
\begin{equation}\label{xi}
 \xi =\partial _t-\lambda^i\partial_x+ \sum_{j\ne i}(\lambda^j-\lambda^i)u^j_x\partial_j+D_x(\lambda^i)u^i_x\partial _{u^i_x}.
\end{equation}
Calculating commutators, we get  $n-1$ vector  fields of the distribution:
 $$\eta_k=\frac{1}{\lambda^k-\lambda^i}[\partial _{u^k_x},\xi]=\partial _k+a_{ik}u^i_x\partial_{u^i_x}.$$
 If there is an extra Riemann invariant of order 0 then, supplemented  with  $\partial _{u^i_x}$, the distribution is involutive,  $\eta_k$ reduces to $\partial_k$ and we have $[\partial_k,\partial _t-\lambda^i\partial_x]=0$ by the Frobenius Theorem. Hence a).

 Suppose that $\lambda^i_j\ne 0$ for some $j\ne i$. Then there is only one Riemann invariant of order 0. We use  $\eta^k$ to kill coefficients of $\partial_j$ and to reduce $\xi$ to
 $$
 \bar{\xi}=\partial_t-\lambda^i\partial_x+\lambda^i_i(u^i_x)^2\partial_{u^i_x}.
 $$

  By the Frobenius theorem, there is a Riemann invariant of order 1 if and only if the distribution $\{\bar{\xi},\partial _{u^k_x},\eta_k,\zeta \}_{k\ne i}$, where
  $$
  \zeta=\frac{1}{\lambda^i_j}[\eta_j,\bar{\xi}]=-\partial_x+b_{ij}(u^i_x)^2\partial_{u^i_x},
  $$
 is involutive.   Computing commutators $[\eta_j,\eta_k]$, $[\eta_j,\zeta]$, we get the first and the second relations of b), comparing $\frac{1}{\lambda^i_k}[\eta_k,\bar\xi]$ with $\zeta$ and computing    $[\eta_k,\frac{1}{\lambda^i_k}[\eta_k,\bar\xi]]$, we obtain the third condition  b). \hfill $\Box$\\
 
  \medskip
\noindent {\bf Remark.} The condition b) of Lemma  \ref{order01} implies that if   $b_{ik}\ne 0$ then the $i$-th
 Lam\'e coefficient, defined as a solutions to the following system of PDEs
\begin{equation}\label{Lame}
\partial_k\ln (H_i)=a_{ik}, \ k \ne i,
\end{equation}
can be chosen as
\begin{equation}\label{explicitLame}
H_i=-\frac{1}{b_{ik}}.
\end{equation}
Note that equations (\ref{Lame}) are compatible due to (\ref{compatibilityCom}), its general solution is defined up to multiplcation by an arbirtary function of $u^i$. Moreover, this system is integrable by quadratures. In fact, let $i=1$ and $u_0=(u^1_0,...,u^n_0)$. Then $\ln (H_1)=\int a_{1n}du^n+\ln (\tilde{H}_1(u^1,u^2,...,u^{n-1}))$ for some $\tilde{H}_1$. For $k\ne 1,n$ we have $\partial_k \ln (H_1)=a_{1k}=\int \partial_k(a_{1n})du^n+\partial_k\ln(\tilde{H}_1)=\int \partial_n(a_{1k})du^n+\partial_k\ln(\tilde{H}_1)=a_{1k}-a_{1k}|_{u^n=u^n_0}+\partial_k\ln(\tilde{H}_1)$. With $\tilde{a}_{1k}(u^1,u^2,...,u^{n-1})=a_{1k}|_{u^n=u^n_0}$ one has $\partial_k\ln(\tilde{H_1})=\tilde{a}_{1k}$ and the claim follows by induction.
 
\section{Darboux integrable diagonal systems of hydrodynamic type}\label{sDi}
System (\ref{hydro}) is  Darboux integrable in order 1 if and only if its characteristic speeds $\lambda^i$ satisfy the conditions b) of lemma (\ref{order01}) for all $i$, these conditions forming a system of differential equations. It turns out that the space of solutions to this system is infinite dimensional. Lets us consider a generic solution and suppose that $ \lambda^i_j\ne 0,$ whenever $ j\ne i$.
\begin{theorem}
Let $\Lambda^i_j$, $\beta_i$, $0\le i,j\le n$ be  germs of real analytic functions of one variables at $u_0^i$. Then there is  (a germ of) system (\ref{hydro}),  Darboux integrable in order 1, such that
$$
\lambda^i_j|_{u^k=u^k_0,\ k\ne i}=\Lambda^i_j,\  \  \  \  b^i_j|_{u^k=u^k_0,\ k\ne i}=\beta_i.
$$
\end{theorem}
{\it Proof:} From semihamiltonian property we have $\forall \ k\ne j\ne i \ne k,$
$$
\lambda^i_{jk} =\frac{(2\lambda^i-\lambda^j  - \lambda^k)\lambda^i_j\lambda^i_k}{(\lambda^i-\lambda^k)}+\frac{(\lambda^k-\lambda^i)\lambda^i_j\lambda^j_k}{(\lambda^k-\lambda^j)(\lambda^j-\lambda^i)}+\frac{(\lambda^j-\lambda^i)\lambda^i_k\lambda^k_j}{(\lambda^j-\lambda^k)(\lambda^k-\lambda^i)}.
$$
For any $i$, the last conditions b) of lemma (\ref{order01}) allow to introduce  functions $B_i=b_{ik}$, satisfying
$$ \partial_kB_i=-a_{ik}B_i$$
 for $k\ne i$. Definition of $b_{ij}$ gives
 $$
 \lambda^i_{ik}= \lambda^i_kB_i-\lambda^i_ia_{ik}
 $$
  for $k\ne i$.
 Invoking the definitions of  $a_{ij}$, we see that  we have a system of partial differential equations for
  functions $\lambda^i, \ \lambda^i_j$, and  $B_i$. By direct calculation one shows that this exterior system is in involution and, by Riquier's  theory \cite{R-10}, the initial conditions specified by the theorem fix a local analytic solution.
\hfill $\Box$\\

Solution of a Darboux integrable semihamiltonian system reduces to solving a system of ordinary differential equations, like for any other Darboux integrable equation or system: find a second Riemann invariant $I^i(x,t,u,u^i_x)$ for each Monge distribution, set
\begin{equation}\label{orbit}
I^i(x,t,u,u^i_x)=\varphi_i(u^i),
\end{equation}
  where $\varphi_i$ are arbitrary,  resolve these equations for $u^i_x$: $u^i_x=P^i(x,t,u,\varphi_i(u^i))$ and integrate the  exterior system
\begin{equation}\label{FrobeniusGeneral}
du^i=P^i(x,t,u,\varphi_i(u^i))(dx+\lambda^i(u)dx).
\end{equation}
\begin{proposition}\label{FrobeniusG}
If all Monge distributions of system (\ref{hydro}) have an extra Riemann invariant  then the exterior system (\ref{FrobeniusGeneral}), where $i=1,...,n,$ is completely integrable (in the sense of Frobenius) for any choice of smooth functions $\varphi_i$.

Moreover, if all but $n$-th  Monge distributions  have an extra Riemann invariant then the exterior system
(\ref{FrobeniusGeneral}), where $i=1,...,n-1,$ supplemented by one equation
$$
du^n=u^n_x(dx+\lambda^n(u)dt),
$$
  admits a Cauchy characteristic for any choice of smooth functions $\varphi_i$.
\end{proposition}
 {\it Proof:} By differentiation of the identity $I^i(x,t,u,P^i(x,t,u,p))=p$ with respect to $x$, $t$, $p$, and  all $u^k$, one finds  partial derivatives of $P^i$.  Now the condition that $I^i$ is invariant along the $i$-th Monge distribution implies that the system (\ref{FrobeniusGeneral}) with $i=1,...,n$ is in involution.

One proves the second claim by checking that the equations
$$
dx+\lambda^ndt=0,\ \ \ du^n_x=u^n_x\left(\lambda^n_nu^n_x+\sum_{k=1}^{n-1}\lambda^n_kP^k\right)dt
$$ define a Cauchy characteristic.
  \hfill $\Box$\\

 \medskip
\noindent {\bf Remark 1.} The second claim of Proposition \ref{FrobeniusG}  follows from the general theory (see \cite{V-28,S-00}): if  $\dim \mathcal{V}$ is odd and $\dim [\mathcal{V},\mathcal{V}]=\dim \mathcal{V}+1$ there is at least one Cauchy characteristic.

  \medskip
\noindent {\bf Example 1: $\bf b_{ik}=0, \ \forall\ k\ne i$.} In this case  one can find the second Riemann invariant by quadratures. The $i$-th Monge distribution is spanned by the following set of vector fields
$$
\{ \partial_x, \ \partial_t+(u^i_x)^2\lambda^i_i\partial_{u^i_x},\ \partial_k+u^i_xa_{ik}\partial_{u^i_x}\}_{k\ne i}.
$$
Then the $i$-th distribution has the invariant $\left(\frac{1}{u^i_x}+\lambda^i_i(u)t \right)H_i(u)$, where  $H_i$ is a Lam\'e coefficient defined by (\ref{Lame}). 
Thus for the system (\ref{FrobeniusGeneral}),  we find $$P^i=\frac{H_i(u)}{\varphi _i(u^i)-t\lambda^i_i(u)H_i(u)}.$$

 \medskip
\noindent {\bf Example 2:  $\bf n=2,$ $\bf b_{12}=b_{21}=0$.} Resolving (\ref{FrobeniusGeneral}) with respect to $dt$ and $dx$, we get:
\begin{eqnarray}
dt=\frac{\varphi_1 du^1}{H_1(\lambda^1-\lambda^2)}+\frac{\varphi_2 du^2}{H_2(\lambda^2-\lambda^1)}+\left( \frac{\lambda^1_1 du^1}{\lambda^2-\lambda^1} +\frac{\lambda^2_2 du^2}{\lambda^1-\lambda^2} \right) t,\\
dx=\frac{\lambda^2\varphi_1 du^1}{H_1(\lambda^2-\lambda^1)}+\frac{\lambda^1 \varphi_2 du^2}{H_2(\lambda^1-\lambda^2)}+\left( \frac{\lambda^2 \lambda^1_1 du^1}{\lambda^1-\lambda^2} +\frac{\lambda^1 \lambda^2_2 du^2}{\lambda^2-\lambda^1} \right) t.
\end{eqnarray}
One checks that the form-valued coefficient of $t$  in the first equation, namely, $\omega=\frac{\lambda^1_1 du^1}{\lambda^2-\lambda^1} +\frac{\lambda^2_2 du^2}{\lambda^1-\lambda^2} $ is closed due to $ b_{12}=b_{21}=0$. Applying a textbook technique of  integration of  linear ODEs, we get
$$
t=W(u) \int \frac{1}{W}\left(\frac{\varphi_1 du^1}{H_1(\lambda^1-\lambda^2)}+\frac{\varphi_2 du^2}{H_2(\lambda^2-\lambda^1)}\right),
$$
where $W=\exp(\int \omega)$. With this $t$, the form for $dx$ does not depend on $x$, therefore the general solution is obtained by quadratures.

 \medskip
\noindent {\bf Remark 2.} Linearly degenerate systems, defined as satisfying $\lambda^i_i=0\ \forall \ i$, are Darboux integrable, provided that they are semihamiltonian,  since for such systems  $b_{ik}=0$. Moreover, the expressions for  $P^i$ can be  rewritten as $P^i=\tilde{\varphi}_i(u^i)H_i(u)$. These formulas appeared in \cite{T-85}, where such system were integrated via (\ref{FrobeniusGeneral}). A different approach, restricted to hamiltonian systems (see \cite{DN-83} for definitions), was given in \cite{P-87}: it applies (\ref{Tsarevform}) and takes advantage of particular properties of conservation laws for linearly degenerate systems.  Independent method  of integration  was described in \cite{F-90}. It has a definite geometric flavor, put in the language of the web theory \cite{BB-38} as follows: on any solution, the Monge characteristic distributions and the  coordinate lines $x=const$, $t=const$ cut an $(n+2)$-web of rank $n$ \cite{F-90}.  For the general case of Darboux integrable systems (\ref{hydro}), no web-geometric interpretation is known. 
 
\section{Darboux integrability of the system for commuting flow}
The above presented results point out that the properties of system (\ref{commuting}) are primal for the Darboux integrability of (\ref{hydro}). 
In this section, we show, in particular, that one is integrable if and only if the other is integrable (in the sense of Darboux).

\begin{theorem}
If system (\ref{hydro}) is Darboux integrable then any commuting system (\ref{hydroMu}) is also Darboux integrable.
\end{theorem}
 {\it Proof:}
The criterion given in lemma \ref{order01} is written in terms of $a_{ij}$, which are the same for both systems.
  \hfill $\Box$\\
One rewrites equations (\ref{commuting})  as an exterior differential system in
$\mathbb{R}^{3n}$ with coordinates $u^i,\mu^i,\mu^i_i=\partial_i\mu^i$:
\begin{equation}\label{extcommuting}
\sigma_i:=\mu^i_idu^i+\sum_{k\neq i}a_{ik}(\mu^k-\mu^i)du^k-d\mu^i=0.
\end{equation}
Solutions to (\ref{commuting}) are integral  manifolds of the system (\ref{extcommuting}) with $du^1\wedge du^2\wedge ... du^n\ne 0$. These $n$-dimensional manifolds are tangent to the distribution $\cal{V}$ with generators
\begin{equation}\label{Vcomm}
\mathcal{V}:=\left\langle  \partial_i+\mu^i_i\partial_{\mu_i}+\sum_{k\neq i}a_{ik}(\mu^k-\mu^i)\partial_{\mu^k}, \partial_{\mu^i_i} \right\rangle_{i=1,...,n}
\end{equation}
Definition \ref{M} of the Monge distributions is applicable to this  $\mathcal{V}$ and we can talk about the Darboux integrability of the system (\ref{commuting}).

Computing the exterior derivatives of  $\sigma^i$ and taking into account relations \ref{compatibilityCom}, one gets the following equations:
\begin{equation}
\Omega_i:=(d\mu^i_i+\sum_{k\ne i}[a_{ik}\mu^i_i+(\mu^i-\mu^k)(a_{ik,i}-a_{ik}a_{ki})]du^k)\wedge du^i=0.
\end{equation}
\begin{theorem}\label{comm01}
Semihamiltonian system (\ref{hydro}) is Darboux integrable in order 0 or 1  if and  only if the corresponding system of commuting flows (\ref{commuting}) is Darboux integrable in order 0 or 1.
\end{theorem}
{\it Proof:} We find the conditions for (\ref{commuting}) to be Darboux integrable. An integral element $$\chi=\sum_j(\alpha_j\partial_j+\beta_j\partial_{\mu_j}+\gamma_j\partial_{\mu^j_j} )$$ is singular if $\omega_k(\chi)=0$, $k=1,...,n$ and the rank of the linear system $$\omega_k(\Xi)=0,\ \ \Omega_k(\chi,\Xi)=0, \ \ k=1,...,n$$ for $\Xi=\sum_i(x_i\partial_i+y_i\partial_{\mu_i}+z_i\partial_{\mu^i_i} )$ drops at $\chi$. It is immediate that it happens if and only if $\alpha_{i}=0$ for some $i$. Therefore the vectors of the $i$-th Monge distribution annihilate the forms
$$
du^i, \ \ d\mu^i_i+\sum_{j\ne i}a_{ij}[(\mu^i-\mu^j)b_{ik}+\mu^i_i]du^k,
$$
and the distribution  is generated by the vectors
\begin{equation}\label{Tj}
T_j:=\partial_j  +\sum_{k\ne j}a_{kj}(\mu^j-\mu^k)\partial_{\mu^k}+\mu^j_j\partial_{\mu^j}-a_{ij}[(\mu^i-\mu^j)b_{ij}+\mu^i_i]\partial_{\mu^i_i},\ \ \ \partial_{\mu^j_j},
\end{equation}
where $j\ne i.$
Computing the commutators $[\partial_{\mu^j_j, T_j}]=\partial_{\mu^j}$, $j\ne i$ one gets  $n-1$ new generators  $\partial_{\mu^j}$ and replaces $T_j$ by
$$
\tilde{T}_j:=\partial_j  +a_{ij}(\mu^j-\mu^i)\partial_{\mu^i}-a_{ij}[(\mu^i-\mu^j)b_{ij}+\mu^i_i]\partial_{\mu^i_i}.
$$
Now $[\partial_{\mu^j},\tilde{T}_j]=a_{ij}[\partial_{\mu^i}+b_{ij}\partial_{\mu^i_i}]$ and
if there is an extra (to $u^i$) Riemann invariant of order 0, then   $a_{ij}=0$, $\forall j\ne i$ and one may choose $\mu^i$ as this invariant.

Suppose that there is no extra invariant of order 0, then at least for one $j\ne i$ holds true $a_{ij}\ne 0$. Moreover, counting the  dimension of the obtained generators of the distribution, one gets $b_{ik}=b_{ij}$  for all $k\ne i $ such that $b_{ik}\ne 0$. Finally, commutators $[\partial_{\mu^i}+b_{ij}\partial_{\mu^i_i},T_k]$  give the last condition of Lemma \ref{order01}.
\hfill $\Box$\\

In contrast to general case of (\ref{hydro}),  extra Riemann  invariants $v^i$ of Darboux integrable system (\ref{commuting}) can be written explicitly: 
\begin{equation}
v^i=H_i(\mu^i_i-b_{ij}\mu^i),
\end{equation}
where $H_i$ is a Lam\'e coefficient, defined by (\ref{Lame}).  
\begin{theorem}\label{commquadr}
Darboux integrable system (\ref{commuting}), possessing Riemann invariants of first order, is integrable in quadratures. 
\end{theorem}
{\it Proof:}  Supplementing system (\ref{commuting}) with $n$ extra equations 
\begin{equation}\label{muii}
\partial_i \mu^i=b_{ij}\mu^i+\frac{\varphi^i(u^i)}{H_i}, 
\end{equation}
where $\varphi^i$ are arbitrary functions of one variable, we obtain a completely integrable system. (Note that for $b_{ij}\ne 0$ equation (\ref{muii}) can be rewritten as $\partial_i \mu^i=b_{ij}(\mu^i-\varphi^i(u^i))$. )

Along a segment parallel to $u^i$-coordinate line the corresponding equations form, in fact, a system of ordinary differential equations, integrable in quadratures: first one integrates a linear ODE (\ref{muii}) then substitute the obtained $\mu^i$ into $\partial_i\mu^k=a_{ik}(\mu^k-\mu^i)$ for $k\ne i$ and gets a system of decoupled linear ODEs, which are also integrable in quadratures.  Let us fix values of $\mu$ at some pont  $u_0$, since any close point $u$ can be connected to $u_0$ by a chain of segments parallel to coordinate lines, one gets a value of $\mu$ at $u$ by integrating a sequence of liner ODEs. 
 \hfill $\Box$\\

Applying  the method of characteristics, we can  reduce integration of  (\ref{commuting}) to ordinary differential equations also in the case when all but one Monge distributions have an extra Riemann invariant.   Let us consider the case of order 1.
\begin{proposition}
If all but $n$-th  Monge distributions of  system (\ref{commuting}) have an extra Riemann invariant  then the exterior system
\begin{equation}\label{FrobeniusSystem}
d\mu^i=\left(b_{ij}\mu^i+\frac{\varphi^i(u^i)}{H_i}\right)du^i+\sum_{k\ne i}a_{ik}(\mu^k-\mu^i)du^k,
\end{equation}
where $i=1,...,n-1$, supplemented by the following equation
\begin{equation}\label{nth}
d\mu^n=\sum_{k=1}^{n-1}a_{nk}(\mu^k-\mu^n)du^k+\mu^n_ndu^n
\end{equation}
admits $n-1$ independent  Cauchy characteristics  for any choice of smooth functions $\varphi^i$. These characteristics can be chosen as 
$$
\partial_j  +\sum_{k\ne j}a_{kj}(\mu^j-\mu^k)\partial_{\mu^k}+\left(b_{ji}\mu^j+\frac{\varphi(u^j)}{H_j}\right)\partial_{\mu^j}-a_{nj}[(\mu^n-\mu^j)b_{nj}+\mu^n_n]\partial_{\mu^n_n}.
$$
\end {proposition}
{\it Proof:} Exterior derivatives of the forms (\ref{FrobeniusSystem}) vanish modulo (\ref{FrobeniusSystem},\ref{nth}). Computing the derivative of (\ref{nth}) and again taking into account (\ref{FrobeniusSystem},\ref{nth}) we get
$$
(d\mu^n_n +\sum_{j=1}^{n-1}a_{nj}[(\mu^n-\mu^j)b_{nj}+\mu^n_n]du^j )\wedge du^n=0.
$$
Therefore any integral element of our system
 satisfying the following 2 equations
$$
du^n=0,\ \ \ d\mu^n_n=-\sum_{j=1}^{n-1}a_{nj}[(\mu^n-\mu^j)b_{nj}+\mu^n_n]du^j
$$
is a Cauchy characteristic.
 \hfill $\Box$\\

\medskip
\noindent {\bf Remark.} Vessiot \cite{V-39} showed that Goursat's classification \cite{G-98} of Darboux integrable scalar equations of second order with two independent variables is governed by Lie groups (see also \cite{S-00} for modern exposition and \cite{V-01} for generalizations). It turns out that {\it Lie vector field systems} play a similar role for Darboux integrable systems of commuting flows (\ref{commuting}). The distribution corresponding to  (\ref{commuting}) is generated by $2n$ vector fields (\ref{Tj}), where now $j$ runs through all the values: $ j=1,...,n.$
Let us rewrite these vector fields in the  Riemann invariants $u^j,v^j$ and $\mu^j$, $ j=1,...,n$.  Then the distribution generators of the system can be chosen as
$$
\tau_j=\partial_j  + \frac{v^j}{H_j}\partial_{\mu^j} +b_{jk}\mu^j\partial_{\mu^j}- \sum_{k\ne j}a_{kj}\mu^k\partial_{\mu^k}+  \sum_{k\ne j}a_{kj}\mu^j\partial_{\mu^k},\ \ \ \partial_{v^j},\ \ \ j=1,...,n,
$$
while excluding the operators $\tau_i,\partial_{v^i}$ gives generators of the $i$-th Monge distribution. All vector fields have a remarkable form: differentiation by one of the Riemann invariants plus a combination of Lie algebra operators of the affine group $\rm{Aff}(\mathbb R^n)$ , i.e. $\partial_{\mu^j},\mu^k\partial_{\mu^j},\ \ j,k=1,...,n$ with coefficients depending only on the Riemann invariants $u,v$. Thus they constitute  {\it Lie vector field system}. 
\section{Laplace transformation  for commuting flows  and  Darboux integrability}

Let us consider system (\ref{commuting}) on its own right, i.e. we do not suppose that $a_{ij}$ are defined via $\lambda$-s by (\ref{aij}).

\begin{lemma}\cite{T-90}
The system (\ref{commuting}) is in involution if
\begin{equation}\label{compatibilityCom}
\partial_j a_{ki}=a_{ki}a_{ij}+a_{kj}a_{ji}-a_{ki}a_{kj}, \ \ \ \ \forall \ j\ne k\ne i\ne j.
\end{equation}
\end{lemma}
 {\it Proof:}  Computing $\partial_i\partial_j\mu^k$ in two different ways and equating the expressions obtained, we get a linear equation  in $\mu^i,\mu^j,\mu^k$. For the system in involution, the coefficients must vanish. Vanishing of the coefficient of $\mu^i$ gives  (\ref{compatibilityCom}). The coefficient of $\mu^k$ gives $\partial_j(a_{ik})=\partial_k(a_{ij})$, which follows from relations (\ref{compatibilityCom}). Finally, the coefficient ob $\mu^j$ gives again (\ref{compatibilityCom}) after  permutation of indexes.
 \hfill $\Box$\\

Let us fix $i$ and suppose that $a_{ik}\ne 0$ for any $k\ne i$. Then we can exclude all $\mu^k$ from system (\ref{commuting}) and obtain a system of compatible PDEs of second order for $\mu^i$ in the following way. Differentiating
\begin{equation}\label{Laplace}
\mu^k=\frac{\partial_k\mu^i}{a_{ik}}+\mu^i
\end{equation}
by $u^m$ with $m\ne i$, $m\ne k$, taking into account (\ref{commuting}), and using subindexes for partial derivatives of $\mu$, we obtain
\begin{equation}\label{Zkim}
 \mu^i_{km}+\left(a_{km}-\frac{\partial_ma_{ik}}{a_{ik}}\right)\mu^i_k+a_{ik}\left( 1-\frac{a_{km}}{a_{im}} \right)\mu^i_m=0.
\end{equation}
Similarly, differentiating (\ref{Laplace}) by $u^i$  we get
\begin{equation}\label{Zkik}
 \mu^i_{ki}+\left(a_{ki}-\frac{\partial_ia_{ik}}{a_{ik}}\right)\mu^i_k+a_{ik}\mu^i_i=0.
\end{equation}
Thus for $U=\mu^i$ we have a system of linear PDEs of the form
\begin{equation}\label{2ndOrder}
U_{ij}+A_{ij}U_i+A_{ji}U_j=0, \ \ \ \ \ i\ne j,
\end{equation}
where again we use the notation $U_i=\partial_iU$, $U_{ij}=\partial_i\partial_jU$.

Conversely, given $i$, we recover  system (\ref{commuting}) from (\ref{2ndOrder}) for
$$
\mu^i=U,\ \ \  \mu^k=\frac{\partial_k U}{A_{ki}+\frac{\partial _i A_{ik}}{A_{ik}}}+U,
$$
with
\begin{equation}\label{aim}
a_{km}=A_{km}+\frac{\partial _mA_{ik}}{A_{ik}}
=A_{im}\left( 1-\frac{A_{mk}}{A_{ik}} \right),
\ \ \mbox{if} \  \ \   i\ne m\ne k\ne i
\end{equation}
and
\begin{equation}\label{aik}
a_{ik}=A_{ik},\ \ \ \  a_{ki}=A_{ki}+\frac{\partial _iA_{ik}}{A_{ik}},\ \  \mbox{if} \  \ i\ne k.
\end{equation}
The two expressions for $a_{km}$ manifest the compatibility conditions for (\ref{2ndOrder}).

A crucial observation is that if we pass from (\ref{2ndOrder}) to (\ref{commuting}) via (\ref{aim},\ref{aik}) for one $i$ and then return from (\ref{commuting}) to another system of type (\ref{2ndOrder}) via (\ref{Laplace}) fixing another index $j\ne i$, we get a Laplace transformation  of (\ref{2ndOrder}) for the chosen pair $i,j$.  Thus, the system for commuting flows (\ref{commuting}) can be viewed as the set of possible Laplace transformations.
But we can as well start with the system (\ref{commuting}). 
\begin{definition}\label{defLT}
Suppose that for a pair of indexes $i,j$ holds true $a_{ik}\ne 0, \ \forall \ k\ne i $ and $a_{jk}-\frac{\partial_k a_{ij}}{a_{ij}}\ne 0, \ \forall \ k\ne j$. Laplace transformation of  system  (\ref{commuting}) for a pair of indexes $i,j$ is obtained by passing to  (\ref{2ndOrder}) for $U=\mu^i$ and returning to (\ref{commuting})  via (\ref{aik},\ref{aim}) with $i=j$:
\begin{equation}\label{LaplaceMu}
\bar{\mu}^j=\mu^i, \ \ \ \bar{\mu}^i=\frac{\mu^i_i}{a_{ji}-\frac{\partial_i a_{ij}}{a_{ij}}}+\mu^i, \ \ \ \bar{\mu}^k=\frac{\mu^i_k}{a_{jk}-\frac{\partial_k a_{ij}}{a_{ij}}}+\mu^i=\frac{a_{ij}\mu^k-a_{kj}\mu^i}{a_{ij}-a_{kj}}, \ \ \ k\ne i\ne j\ne k.
\end{equation}
\end{definition}
Applying formulas (\ref{Zkim},\ref{Zkik},\ref{aim},\ref{aik}), one computes all the coefficients $\bar{a}_{km}$ of the transformed system.

\medskip

\noindent {\bf Remark.} Since $\lambda$ is a solution to system (\ref{commuting}) the Laplace transformation acts also on generic solutions to system (\ref{hydro}): one applies formulas (\ref{LaplaceMu}) to $\lambda$-s and invokes the implicit form (\ref{Tsarevform}). Note that this action is not a point transformation.

\medskip

There are two ways for Laplace transformation introduced by Definition \ref{defLT} to fail. First, by passing from (\ref{commuting}) to (\ref{2ndOrder}) the obstacle is vanishing of $a_{ik}$ for some or all $k$. Second, with $a_{ji}-\frac{\partial_i a_{ij}}{a_{ij}}\ne 0$, returning to  (\ref{commuting}) is impossible if $a_{jk}-\frac{\partial_k a_{ij}}{a_{ij}}=0$ for some or all $k$.

Consider the first obstacle. Suppose that the set of indexes ${\cal{I}}_i=\{1,2,...,n\}\setminus \{i\}$ is a union of two disjunctive subsets ${\cal{I}}_i={\cal{Z}}\cup {\cal{N}}$.   For $\alpha\in  {\cal{Z}}$ holds $a_{i\alpha}= 0$  and for $\beta\in  {\cal{N}}$ holds $a_{i\beta}\ne 0$. 

\begin{proposition}\label{subsystem}
The functions $\mu^i$ and $\mu^{\beta}, \ \beta\in  {\cal{N}}$ satisfy a subsystem of the form (\ref{commuting}), namely $\forall \alpha\in  {\cal{Z}}$   and $\forall \beta,\bar{\beta}\in  {\cal{N}}, \ \beta\ne \bar{\beta}$ holds
$$
 \begin{array}{l}
  \partial_{\alpha} \mu^i=0,\\
 \partial_{\alpha} \mu^{\beta}=0,\\
 \partial_{\beta} \mu^i=a_{i\beta}(\mu^{\beta}-\mu^i), \ \ \mbox{where}\  \partial_{\alpha}a_{i\beta}=0,\\
  \partial_i \mu^{\beta}=a_{\beta i}(\mu^i-\mu^{\beta}), \ \ \mbox{where}\  \partial_{\alpha}a_{\beta i}=0,\\
   \partial_{\beta} \mu^{\bar{\beta}}=a_{\bar{\beta}\beta}(\mu^{\beta}-\mu^{\bar{\beta}}), \ \ \mbox{where}\  \partial_{\alpha}a_{\bar{\beta}\beta}=0.
 \end{array}
 $$
\end{proposition}
 {\it Proof:}
Compatibility conditions  (\ref{compatibilityCom}) (ou direct comparing of the mixed derivatives of $\mu^i$ by $u_{\beta}$ e $u_{\alpha}$ found from (\ref{commuting})) and the inequality $a_{i\beta}\ne 0$ imply $\partial_{\alpha}a_{i\beta}=0$ and $a_{\beta \alpha}=0$. Therefore $\partial_{\alpha} \mu^{\beta}=0$. The semihamiltonian property gives  $\partial_{\alpha}a_{\beta i}=0$ and $\partial_{\alpha}a_{\bar{\beta}\beta}=0$.
  \hfill $\Box$\\
\begin{lemma}\label{compatibilityA}
 Coefficients $A_{ki}$ of equations (\ref{2ndOrder}) satisfy the following compatibility relations, similar to (\ref{compatibilityCom})  
$$
\partial_j A_{ki}=A_{ki}A_{kj}-A_{ki}A_{ij}-A_{ji}A_{kj}, \ \ \ \ \forall \ j\ne k\ne i\ne j.
$$
In particular, $\partial_j A_{ki}=\partial_i A_{kj}.$
\end{lemma}
 {\it Proof:} 
We  differentiate $U_{ij}$, expressed in terms of $U_i,U_j$,  by $u_k$   and $U_{ik}$,  expressed in terms of $U_i,U_k$,  by $u_j$, equate the results, and isolate the coefficient of $U_k$. \hfill $\Box$\\

Now let us study the second obstacle. Suppose that for a pair of indexes $i,j$ holds true $a_{ik}\ne 0,\ \forall  k\ne i$ and   $b_{ij}\ne 0$. 
Again, we represent the set of indexes ${\cal{I}}_j=\{1,2,...,n\}\setminus \{j\}$ as  a union of two disjunctive subsets ${\cal{I}}={\bar{\cal{Z}}}\cup \bar{{\cal{N}}}$ so that  for $\alpha\in  \bar{{\cal{Z}}}$ holds $a_{j\alpha}-\frac{\partial_{\alpha} a_{ij}}{a_{ij}}= 0$  and for $\beta\in \bar{ {\cal{N}}}$ holds $a_{j\beta}-\frac{\partial_{\beta} a_{ij}}{a_{ij}}\ne 0$. Note that $i\in \bar{ {\cal{N}}}$. 

\begin{proposition}
There is a system of the form (\ref{commuting}) for $\bar{\mu}^j$ and $\bar{\mu}^{\beta}, \ \beta \in \bar{ {\cal{N}}}$

$$
 \begin{array}{l}
 \partial_{\beta} \bar{\mu}^j=\bar{a}_{j\beta}(\bar{\mu}^{\beta}-\bar{\mu}^j), \ \ \mbox{where}\  \partial_{\alpha}\bar{a}_{j\beta}=0,\\
  \partial_j \bar{\mu}^{\beta}=\bar{a}_{\beta j}(\bar{\mu}^j-\bar{\mu}^{\beta}), \ \ \mbox{where}\  \partial_{\alpha}\bar{a}_{\beta j}=0,\\
   \partial_{\beta} \bar{\mu}^{\bar{\beta}}=\bar{a}_{\bar{\beta}\beta}(\bar{\mu}^{\beta}-\bar{\mu}^{\bar{\beta}}), \ \ \mbox{where}\  \partial_{\alpha}\bar{a}_{\bar{\beta}\beta}=0, \ \bar{\beta}\in \bar{N}.
 \end{array}
 $$
such that excluding $\bar{\mu}^{\beta}, \ \beta\ne j$ gives the subset of equations (\ref{2ndOrder}) for $U=\bar{\mu}^j$ with low indexes from $\bar{ {\cal{N}}}$.  
\end{proposition}
{\it Proof:} Excluding $\bar{\mu}^k, \ k\ne j$ from the system to be constructed via $$
\bar{\mu}^k=\frac{\partial_k\bar{\mu}^j}{\bar{a}_{jk}}+\bar{\mu}^j
$$
 is possible only if $\bar{a}_{jk}\ne 0$. The resulting system of the form (\ref{2ndOrder}) for $U=\bar{\mu}^j$ contains the following equations  
$$
 U_{kj}+\left(\bar{a}_{kj}-\frac{\partial_j\bar{a}_{jk}}{\bar{a}_{jk}}\right)U_k+\bar{a}_{jk}U_j=0
$$
(apply (\ref{Zkik}) for $i$ replaced by $j$). If these equations coincide with the corresponding equations obtained by excluding $\mu^i$ from (\ref{commuting}) then $\bar{a}_{jk}=a_{jk}-\frac{\partial_k a_{ij}}{a_{ij}}$. Therefore we can not recover a system of the form (\ref{commuting}) for all indexes but only for  $ \beta\in \bar{ {\cal{N}}}$ (and therefore for $j$). 

By (\ref{aim},\ref{aik},\ref{Laplace}), we set 
$$
\bar{a}_{j\beta}=A_{j\beta},\ \ \ \  
\bar{a}_{\beta j}=\frac{\partial_j A_{j\beta}}{A_{j\beta}}+A_{\beta j},\ \ \ \ 
\bar{a}_{\bar{\beta} \beta}=\frac{\partial_{\beta}A_{j\bar{\beta}}}{A_{j\bar{\beta}}}+A_{\bar{\beta}\beta},\ \ \ \ 
\bar{\mu}^{\beta}=\frac{\partial_{\beta}\bar{\mu}^j}{A_{j\beta}}+\bar{\mu}^j. 
$$ 
Then $\partial_{\beta} \bar{\mu}^j=\bar{a}_{j\beta}(\bar{\mu}^{\beta}-\bar{\mu}^j)$ is true by definition of $\bar{\mu}^{\beta}$ and $\bar{a}_{j\beta}$. Substituting the above defined   $\bar{a}_{\beta j}$ and $\bar{\mu}^{\beta}$ into  
$\partial_j \bar{\mu}^{\beta}=\bar{a}_{\beta j}(\bar{\mu}^j-\bar{\mu}^{\beta})$ 
we obtain equations (\ref{2ndOrder}) with subindexes $j,\beta$ for $U=\bar{\mu}^j$. Finally, substituting  
$\bar{\mu}^{\bar{\beta}}$ and $\bar{a}_{\bar{\beta} \beta}$ into $ \partial_{\beta} \bar{\mu}^{\bar{\beta}}=\bar{a}_{\bar{\beta}\beta}(\bar{\mu}^{\beta}-\bar{\mu}^{\bar{\beta}})$ we see that the coefficient of $\partial_{\beta}\bar{\mu}^j$ is equal to $A_{\beta \bar{\beta}}$ do to Lemma \ref{compatibilityA}, thus obtaining equations (\ref{2ndOrder}) with subindexes $\beta,\bar{\beta}$ for $U=\bar{\mu}^j$.

By Lemma (\ref{compatibilityA}), we have $\partial_{\alpha}\bar{a}_{j\beta}=\partial_{\alpha}A_{j\beta}=\partial_{\beta}A_{j\alpha}=0$ since $A_{j\alpha}=0$. 
Applying Lemma (\ref{compatibilityA}) to $\partial_{\beta}A_{j\alpha}$ and taking into account that $A_{j\alpha}=0$, $A_{j\beta}\ne 0$ we get $A_{\beta\alpha}=0.$  
Now substituting the expression for $\bar{a}_{\beta j}$ into  $\partial_{\alpha}\bar{a}_{\beta j}$,  again by  
Lemma  (\ref{compatibilityA}) and due to $A_{\beta\alpha}=0$, we obtain   $\partial_{\alpha}\bar{a}_{\beta j}=0$. 
Similar, $A_{\bar{\beta}\alpha}=0$, $A_{j\alpha}=0$ imply   $\partial_{\alpha}\bar{a}_{\bar{\beta}\beta}=0$.
 \hfill $\Box$\\

\begin{definition}\label{defRLT}
Suppose that for a pair of indexes $i,j$ holds true $a_{ij}\ne 0, \ b_{ij}\ne 0$. Reduced Laplace transformation of  system  (\ref{commuting}) for a pair of indexes $i,j$ is obtained by passing  from the subsystem for  $\mu^{\beta}, \ \beta\in  {\cal{N}}$ to  (\ref{2ndOrder}) for $U=\mu^i$ and recovering the subsystem of the form  (\ref{commuting}) for $\bar{\mu}^j=\mu^i$, $\bar{\mu}^{\beta}, \ \beta\in  \bar{{\cal{N}}}\cap {\cal{N}} $:
\begin{equation}\label{RLaplaceMu}
\bar{\mu}^j=\mu^i, \ \ \ \bar{\mu}^i=\frac{\mu^i_i}{a_{ji}-\frac{\partial_i a_{ij}}{a_{ij}}}+\mu^i, \ \ \ \bar{\mu}^{\beta}=\frac{\mu^i_{\beta}}{a_{j\beta}-\frac{\partial_{\beta} a_{ij}}{a_{ij}}}+\mu^i=\frac{a_{ij}\mu^{\beta}-a_{\beta j}\mu^i}{a_{ij}-a_{\beta j}}, \ \ \ \beta\ne i\ne j\ne \beta.
\end{equation}
\end{definition} 
Thus the reduced Laplace transformation is passing from (\ref{commuting}) to a system of the same form with the smaller number of equations and field variables.  
\begin{theorem}\label{terminates1}
Suppose that for a pair of indexes $i,j$ holds true $a_{ij}\ne 0, \ b_{ij}\ne 0$.  Then its $i$-th Monge distribution of  semihamiltonian system (\ref{hydro}) has an extra Riemann invariant of order 1 if and only if all the coefficients $\bar{a}_{i\beta}$  of its reduced Laplace transformation   vanish:  $\bar{a}_{i\beta}=0$.
\end{theorem}
 {\it Proof:} Due to Proposition \ref{subsystem}, we can suppose that the $a_{ik}\ne 0\ \forall k\ne i$. Two systems (\ref{commuting}) are related by the reduced $ij$- Laplace transformation if  equations (\ref{2ndOrder}) for $U=\mu^i$ and $U=\bar{\mu}^j$ are the same. Thus, by (\ref{aim},\ref{aik}) we have for $m\ne j$:
 $$
 \bar{a}_{im}=A_{im}+\frac{\partial _mA_{ji}}{A_{ji}}=a_{im}+\frac{\partial _m(a_{ji}-\frac{\partial_i a_{ij}}{a_{ij}})}{a_{ji}-\frac{\partial_i a_{ij}}{a_{ij}}}.
 $$
On the other hand, using another expression for $\bar{a}_{im}$, we also have
 $$
 \bar{a}_{im}=A_{jm}\left(1-\frac{A_{mi}}{A_{ji}} \right)=a_{im}\left(1-\frac{a_{mj}}{a_{ij}} \right)\left(1-\frac{a_{mi}-\frac{\partial_i a_{im}}{a_{im}}}{a_{ji}-\frac{\partial_i a_{ij}}{a_{ij}}} \right).
 $$
 Similarly,
 $$
  \bar{a}_{ij}=A_{ij}+\frac{\partial _jA_{ji}}{A_{ji}}=a_{ij}+\frac{\partial _j(a_{ji}-\frac{\partial_i a_{ij}}{a_{ij}})}{a_{ji}-\frac{\partial_i a_{ij}}{a_{ij}}}.
 $$
 Thus, vanishing of $ \bar{a}_{im}$ for $m\ne j$ gives
 $$
 a_{im}+\frac{\partial_m b_{ij}}{b_{ij}}=0,\ \ \ \ b_{im}=b_{ij},
 $$
 vanishing of $ \bar{a}_{ij}$ yields
 $$
 a_{ij}+\frac{\partial_j b_{ij}}{b_{ij}}=0,
 $$
 and the theorem follows from
 lemma \ref{order01}. \hfill $\Box$\\

 The above Theorem and Lemma \ref{order01} motivates the following definition.
 \begin{definition}\label{defterminates}
 We say that  $i$-th sequence of reduced Laplace transformation for a semihamiltonian system $S_0$  terminates after $r$ steps if there is a sequence $j_1,j_2,...,j_r\in \mathbb{N}$ of indexes and a sequence $S_1,S_2,...,S_r$ of  semihamiltonian systems such that $S_q$ is a reduced $ij_q$-Laplace transformation of $S_{q-1}$ and either $a_{im}^{(r)}=0, \ \forall m\ne i$ for the system $S_r$   or $b_{im}^{(r-1)}=0, \ \forall m\ne i$ for the system $S_{r-1}$. 
 \end{definition}

Later we show that 
the $i$-th Monge distribution of a semihamiltonian system (\ref{hydro}) has an extra Riemann invariant of order $r$ if and only if any its $i$-th sequence of reduced Laplace transformation terminates after $r$ steps.

 For the case $n=2$,  the hodograph transform reduces a two-component system (\ref{hydro}) to a single equation of form (\ref{2ndOrder}), for which the result is well-known.
 
 \section{Darboux integrability for prolongations}
in this section, we extend the obtained results to the prolongations of systems (\ref{hydro}) and (\ref{commuting}).
The key  to Darboux integrability of prolonged systems  is  the property of being semihamiltonian.

\begin{theorem}\label{Darbouxissemihamiltonian}
 If the $i$-th Monge distribution of a strictly hyperbolic system  (\ref{hydro}) has an extra Riemann invariant of order $r$  then  $\partial_ja_{ik}=\partial_ka_{ij}$ $\forall \ k\ne j\ne i\ne k.$ In particular, Darboux integrable  strictly hyperbolic  system (\ref{hydro})  is semihamiltonian.
 \end{theorem}
 {\it Proof:} Let us fix $i$, the order $r>1$ and introduce the notation $u^k_m:=D_x^m(u^k)$. The $i$-th Monge distribution of the $r$-th prolongation of (\ref{hydro}) is generated by the following vector fields:
 \begin{equation}\label{xir}
 \xi =\partial _t-\lambda^i\partial_x+ \sum_{k\ne i}\sum_{m=0}^{r-1}[D_x^m(\lambda ^ku^k_1)-\lambda ^iu^k_{m+1}]\partial _{u^k_m}+\sum_{m=1}^{r}[D_x^m(\lambda ^iu^i_1)-\lambda ^iu^i_{m+1}]\partial _{u^i_m},
\end{equation} which is the properly truncated operator $D_t-\lambda^iD_x$, restricted to the system, and by   $n-1$ differentiations
 $\eta^{\alpha}_r:=\partial_{u^{\alpha}_r}$ with ${\alpha}\ne i$. (In what follows in this proof, we suppose ${\alpha}\ne i$.)
  To compute or to guess a nice set of generators for Monge distribution of the prolonged system does not seem very promising since the order $r$ may be arbitrary large. Our scheme is to keep track of some "important" variables in recursive determination of the vector fields
  $$
  \eta^{\alpha}_{r-(s+1)}:=\frac{1}{\lambda^{\alpha}-\lambda^i}[\eta^{\alpha}_{r-s},\xi]
  $$
  and finally extract the desired claim from the commutator
  $$
  [\eta^{\beta}_{r-1},\eta^{\alpha}_0]=[\partial_{u^{\beta}_{r-1}}+a_{i\beta}u^i_1\partial_{u^i_r},\eta^{\alpha}_0].
  $$
 Thus, we have to understand how the variables $u^i_r$ and $u^k_{r-1}$ with $k\ne i$ enter in the coefficients of $\eta^{\alpha}_0$.

 We will use the following commutation relation, easily proved by induction:
 $$
 \partial_{u^k_p}D_x^m=\sum_{s=0}^mC_m^sD_x^{m-s}\partial_{u^k_{p-s}},
 $$
where $C_m^s$ are the binomial coefficients. 

Applying the above relation for $\alpha\ne i$ and $p\ge 1$, one gets
\begin{equation}\label{daxi}
\begin{array}{l}
\frac{1}{\lambda^{\alpha}-\lambda^i}[\partial_{u^{\alpha}_p},\xi]=\partial_{u^{\alpha}_{p-1}}+a_{i\alpha}u^i_1\partial_{u^i_p}
+\sum_{m=p+1}^r\frac{C_m^pD_x^{m-p}(\lambda^i_{\alpha}u^i_1)}{\lambda^{\alpha}-\lambda^i}\partial_{u^i_m}+\\
\\
+\sum_{m=p}^{r-1}\frac{C_m^{p-1}D_x^{m-p+1}(\lambda^{\alpha})}{\lambda^{\alpha}-\lambda^i}\partial_{u^{\alpha}_m}
+\sum_{s\ne i}\sum_{m=p}^{r-1}\frac{C_m^pD_x^{m-p}(\lambda^s_{\alpha}u^s_1)}{\lambda^{\alpha}-\lambda^i}\partial_{u^s_m}
\end{array}
\end{equation}
and similarly for $p\ge 2$
\begin{eqnarray*}
[\partial_{u^i_p},\xi]=\sum_{m=p}^rC_m^pD_x^{m-p}(\lambda^i_iu^i_1)\partial_{u^i_m}+
\sum_{m=p}^{r}C_m^{p-1}D_x^{m-p+1}(\lambda^i)\partial_{u^i_m}
+\sum_{s\ne i}\sum_{m=p}^{r-1}C_m^pD_x^{m-p}(\lambda^s_iu^s_1)\partial_{u^s_m}.
\end{eqnarray*}
First we note that the operator $\xi$ does not raise the order of the derivatives $u^i_m$ whereas it increases the order of  $u^{\alpha}_m$ by one.
From  the above commutator formulae one easily sees that  $u^i_r$ appears only in the last step of the recursion generation of $\eta^{\alpha}_0$ and
$$
\eta^{\alpha}_0=\partial_{u^{\alpha}}+a_{i\alpha}u^i_1\partial_{u^i_1}+ra_{i\alpha}u^i_r\partial_{u^i_r}+...
$$
The following observations allow to trace $u^{\beta}_{r-1}$.\\
\begin{itemize}
\item The operator $\eta^{\alpha}_{r-s}$ has the form  $\eta^{\alpha}_{r-s}=\partial_{u^{\alpha}_{r-s}}+a_{i\alpha}u^i_1\partial_{u^i_{r-s+1}}+...$, where omitted are the derivations by $u^{\alpha}_p$ with $p>r-s$ and by $u^i_p$ with $p>r-s+1$.
\item The operator $\xi$ raises the order $s-1$ of $u^{\alpha}_{s-1}$, present in $\eta^{\alpha}_{r-s}$, to order $r-1$ of $u^{\alpha}_{r-1}$ in $\eta^{\alpha}_0$.
\item The derivative $u^{\alpha}_{s-1}$ of order  $s-1$ appears in $\eta^{\alpha}_{r-s}$ by exactly three mechanisms. The first is action of $\xi$ on $u^{\alpha}_{s-2}$, present in $\eta^{\alpha}_{r-s-1}$.
\item The second is generating by $[\partial_{u^{\alpha}_{r-s-1}},\xi]$.
\item The third is action of $\xi$ on $u^i_{s-1}$ present in $\eta^{\alpha}_{r-s-1}$.
\end{itemize}
Let us compute $\eta^{\alpha}_{r-2}$ keeping only the derivations', generating the "important" variables, and  the terms containing $u^{s}_1$ with $s\ne i$ and $u^i_2$. Using the above presented commuting relations we get
\begin{eqnarray*}
\eta^{\alpha}_{r-2}=\partial_{u^{\alpha}_{r-2}}+a_{i\alpha}u^i_1\partial_{u^i_{r-1}}+
C^{r-2}_{r-1}\sum_{k\ne i}u^k_1\frac{  \lambda^{\alpha}_k}{\lambda^{\alpha}-\lambda^i}\partial_{u^{\alpha}_{r-1}}+
C^{r-1}_{r-1}\sum_{k\ne i}u^k_1\frac{  \lambda^k_{\alpha}}{\lambda^{\alpha}-\lambda^i}\partial_{u^k_{r-1}}+\\
\\
+u^i_1\sum_{k\ne i}u^k_1\left[ C^{r-1}_r\frac{  \lambda^i_{\alpha k}}{\lambda^{\alpha}-\lambda^i}+(C^{r-1}_r-1)\frac{ a_{i\alpha} \lambda^i_k}{\lambda^{\alpha}-\lambda^i} +\frac{\lambda^k-\lambda^i}{\lambda^i-\lambda^{\alpha}} \partial_ka_{i \alpha} \right]\partial_{u^i_r}+C^{r-1}_ra_{i\alpha}u^i_2\partial_{u^i_r}.
\end{eqnarray*}
Taking into account the observations, we deduce that the operator $\eta^{\alpha}_{r-s}$ has the form
\begin{eqnarray*}
\eta^{\alpha}_{r-s}=\partial_{u^{\alpha}_{r-s}}+a_{i\alpha}u^i_1\partial_{u^i_{r-s+1}}+\left(U^{\alpha}_{r-s}u^i_s+ u^i_1\sum_{k\ne i}X^k_{r-s}u^k_{s-1}\right)
\partial_{u^i_r}+\\
\\
+\sum_{k\ne i}Y^k_{r-s}u^k_{s-1}\partial_{u^{\alpha}_{r-1}}+
\sum_{k\ne i}Z^k_{r-s}u^k_{s-1}\partial_{u^k_{r-1}}+...,
\end{eqnarray*}
Where $X^{\alpha}_{r-s}, Y^{\alpha}_{r-s}, Z^{\alpha}_{r-s}, U^{\alpha}_{r-s}$ do depend on derivatives.  Thus we have
\begin{eqnarray*}
\eta^{\alpha}_{r-(s+1)}=\frac{1}{\lambda^{\alpha}-\lambda^i}[\eta^{\alpha}_{r-s},\xi]=
\frac{1}{\lambda^{\alpha}-\lambda^i}[\partial_{u^{\alpha}_{r-s}},\xi]+
\frac{a_{i\alpha}u^i_1}{\lambda^{\alpha}-\lambda^i}[\partial_{u^i_{r-s+1}},\xi]-\frac{1}{\lambda^{\alpha}-\lambda^i}U^{\alpha}_{r-s}\xi(u^i_s)\partial_{u^i_r}+\\
\\
-\frac{u^i_1}{\lambda^{\alpha}-\lambda^i}\sum_{k\ne i}X^k_{r-s}\xi(u^k_{s-1})\partial_{u^i_r}
-\frac{1}{\lambda^{\alpha}-\lambda^i}\sum_{k\ne i}Y^k_{r-s}\xi(u^k_{s-1})\partial_{u^{\alpha}_{r-1}}
-\frac{1}{\lambda^{\alpha}-\lambda^i}\sum_{k\ne i}Z^k_{r-s}\xi(u^k_{s-1})\partial_{u^k_{r-1}}+...
\end{eqnarray*}  
The wanted contribution of the term $\frac{1}{\lambda^{\alpha}-\lambda^i}[\partial_{u^{\alpha}_{r-s}},\xi]$ is 
\begin{eqnarray*}
\partial_{u^{\alpha}_{r-s-1}}+a_{i\alpha}u^i_1\partial_{u^i_{r-s}}+
\frac{C_r^{r-s}u^i_1}{\lambda^{\alpha}-\lambda^i}\sum_{k\ne i}\lambda^i_{\alpha k}u^k_s\partial_{u^i_r}+
\frac{C_{r-1}^{r-s-1}}{\lambda^{\alpha}-\lambda^i}\sum_{k\ne i}\lambda^{\alpha}_k u^k_s\partial_{u^{\alpha}_{r-1}}+
\frac{C_{r-1}^{r-s}}{\lambda^{\alpha}-\lambda^i}\sum_{k\ne i}\lambda^k_{\alpha}u^k_s\partial_{u^{k}_{r-1}}, 
\end{eqnarray*}  
the one from $[\partial_{u^i_{r-s+1}},\xi]$ is 
$$
C^{r-s}_r\sum_{k\ne i}\lambda^i_ku^k_s\partial_{u^i_r},
$$
the one from $\xi(u^k_{s-1})$ is
$$
(\lambda^k-\lambda^i)u^k_s,
$$
and the one from 
$
-\frac{1}{\lambda^{\alpha}-\lambda^i}U^{\alpha}_{r-s}\xi(u^i_s)\partial_{u^i_r}
$
is 
$$
-\frac{C^{r-s+1}_ra_{i\alpha }u^i_1}{\lambda^{\alpha}-\lambda^i}\sum_{k\ne i}\lambda^i_ku^k_s\partial_{u^i_r}
$$
since 
$U^{\alpha}_{r-s}=C^{r-s+1}_ra_{i\alpha }.$
Finally,  taking into account that 
$$
\xi(u^k_{s-1})=(\lambda^k-\lambda^i)u^k_s+...
$$
and collecting similar terms, we get the following recursive equations for 
$X^{\alpha}_{m}, Y^{\alpha}_{m}, Z^{\alpha}_{m}$: 
$$
X^k_{r-(s+1)}=\frac{\lambda^k-\lambda^i}{\lambda^i-\lambda^{\alpha}}X^k_{r-s}+C^{r-s}_r\frac{\lambda^i_{\alpha k}}{\lambda^{\alpha}-\lambda^i}+
C^{r-s}_r\frac{a_{i\alpha}\lambda^i_k}{\lambda^{\alpha}-\lambda^i}-C^{r-s+1}_r\frac{a_{i\alpha}\lambda^i_k}{\lambda^{\alpha}-\lambda^i},\\
$$
$$
Y^k_{r-(s+1)}=\frac{\lambda^k-\lambda^i}{\lambda^i-\lambda^{\alpha}}Y^k_{r-s}+C^{r-s-1}_{r-1}\frac{\lambda^{\alpha}_k}{\lambda^{\alpha}-\lambda^i},
$$ 
$$
Z^k_{r-(s+1)}=\frac{\lambda^k-\lambda^i}{\lambda^i-\lambda^{\alpha}}Z^k_{r-s}+C^{r-s}_{r-1}\frac{\lambda^k_{\alpha}}{\lambda^{\alpha}-\lambda^i}.
$$ 
From the operator $\eta^{\alpha}_{r-2}$, we extract the initial values  
\begin{eqnarray*}
X^k_{r-2}=C^{r-1}_r\left(\frac{  \lambda^i_{\alpha k}}{\lambda^{\alpha}-\lambda^i}+\frac{ a_{i\alpha} \lambda^i_k}{\lambda^{\alpha}-\lambda^i}\right)-\frac{ a_{i\alpha} \lambda^i_k}{\lambda^{\alpha}-\lambda^i} +\frac{\lambda^k-\lambda^i}{\lambda^i-\lambda^{\alpha}} \partial_ka_{i \alpha},\\
\\
Y^k_{r-2}=C^{r-2}_{r-1}\frac{\lambda^{\alpha}_k}{\lambda^{\alpha}-\lambda^i}, \ \ \ \ \ \ 
Z^k_{r-2}=C^{r-1}_{r-1}\frac{\lambda^k_{\alpha}}{\lambda^{\alpha}-\lambda^i}.
\end{eqnarray*}
It is not difficult to "integrate" the above equations and obtain   
\begin{eqnarray*}
X^k_0=\left(\frac{  \lambda^i_{\alpha k}}{\lambda^{\alpha}-\lambda^i}+\frac{ a_{i\alpha} \lambda^i_k}{\lambda^{\alpha}-\lambda^i}\right)\frac{(z_k+1)^r-1-z_k^r}{z_k}
-\frac{ a_{i\alpha} \lambda^i_k}{\lambda^{\alpha}-\lambda^i}\frac{(z_k+1)^r-1-rz_k}{z_k^2}
 +z_k^{r-1}\partial_ka_{i \alpha},\\
\\
Y^k_0=\frac{\lambda^{\alpha}_k}{\lambda^{\alpha}-\lambda^i}[(z_k+1)^{r-1}-z_k^{r-1}], \ \ \ \ \ \
Z^k_0=\frac{\lambda^k_{\alpha}}{\lambda^{\alpha}-\lambda^i}\frac{(z_k+1)^{r-1}-1}{z_k}, \ \ \mbox{where}\ \ z_k=\frac{\lambda^k-\lambda^i}{\lambda^i-\lambda^{\alpha}}. 
\end{eqnarray*}
Now observe that if the vector field $[\eta^{\beta}_{r-1},\eta^{\alpha}_0]$ is in the distribution then as well is $[\eta^{\beta}_{r-1},\bar{\eta}^{\alpha}_0]$, where
\begin{eqnarray*}
\bar{\eta}^{\alpha}_0=\eta^{\alpha}_0-\sum_{k\ne i}u^k_{r-1}Y^k_0\eta^{\alpha}_{r-1}-\sum_{k\ne i}u^k_{r-1}Z^k_0\eta^k_{r-1}=\\
\\
\partial_{u^{\alpha}}+a_{i\alpha}u^i_1\partial_{u^i_1}+\left(ra_{i\alpha}u^i_r + u^i_1\sum_{k\ne i}u^k_{r-1}[X^k_0-a_{i\alpha}Y^k_0-a_{ik}Z^k_0]   \right)\partial_{u^i_r}+...
\end{eqnarray*}
The commutator $[\eta^{\beta}_{r-1},\bar{\eta}^{\alpha}_0]=E\partial_{u^i_r}$ with 
$$
E=X^{\beta}_0-a_{i\alpha}Y^{\beta}_0-a_{i\beta}Z^{\beta}_0+a_{i\alpha}a_{i\beta}(r-1)-\partial_{\alpha}a_{i \beta}
$$
must vanish since the distribution is supposed to have an extra Riemann invariant of order $r$. 
Substituting into the equation $E=0$ the found expressions for $X^{\beta}_{0}, Y^{\beta}_{0}, Z^{\beta}_{0}$ we get
\begin{eqnarray*}
\frac{(z_{\beta}+1)^r}{z_{\beta}}\left( \frac{  \lambda^i_{\alpha \beta}}{\lambda^{\alpha}-\lambda^i}+\frac{ a_{i\alpha} \lambda^i_{\beta}}{\lambda^{\alpha}-\lambda^i}
 -\frac{ a_{i\alpha} \lambda^i_{\beta}}{z_{\beta}(\lambda^{\alpha}-\lambda^i)}
 -\frac{z_{\beta} a_{i\alpha} \lambda^{\alpha}_{\beta}}{(z_{\beta}+1)(\lambda^{\alpha}-\lambda^i)}
 -\frac{a_{i\beta} \lambda_{\alpha}^{\beta}}{(z_{\beta}+1)(\lambda^{\alpha}-\lambda^i)}  \right)+\\
 \\
 z_{\beta}^{r-1}\left( 
 -\frac{  \lambda^i_{\alpha \beta}}{\lambda^{\alpha}-\lambda^i}-\frac{ a_{i\alpha} \lambda^i_{\beta}}{\lambda^{\alpha}-\lambda^i}
 +\partial_{\beta}a_{i\alpha}+\frac{a_{i\alpha}\lambda^{\alpha}_{\beta}}{\lambda^{\alpha}-\lambda^i}
     \right)+
 r\left( \frac{a_{i\alpha}\lambda^i_{\beta}}{z_{\beta}(\lambda^{\alpha}-\lambda^i)} +a_{i\alpha}a_{i\beta} \right)+\\
 \\
 \frac{1}{z_{\beta}} \left( 
 -\frac{  \lambda^i_{\alpha \beta}}{\lambda^{\alpha}-\lambda^i}-\frac{ a_{i\alpha} \lambda^i_{\beta}}{\lambda^{\alpha}-\lambda^i}+
 \frac{ a_{i\alpha} \lambda^i_{\beta}}{z_{\beta}(\lambda^{\alpha}-\lambda^i)}+
 \frac{a_{i\beta} \lambda_{\alpha}^{\beta}}{\lambda^{\alpha}-\lambda^i} - z_{\beta}a_{i\alpha}a_{i\beta} -z_{\beta} \partial_{\alpha}a_{i \beta}
 \right) =0.  
\end{eqnarray*}
Substituting the formulae for $a_{ij}$ and $z_{\beta}$, one checks by direct computation that the coefficients of $z_{\beta}^{r-1}$, $r$, and $\frac{1}{z_{\beta}}$ vanish identically while the vanishing of the coefficient by $\frac{(z_{\beta}+1)^r}{z_{\beta}}$ is equivalent to the semihamiltonian property 
 $\partial_{\alpha}a_{i\beta}=\partial_{\beta}a_{i\alpha }$. 
 \hfill $\Box$
\begin{proposition}
 Suppose that a semihamiltonian system (\ref{hydro}) contains a subsystem, described in Proposition \ref{subsystem}, and has an extra Riemann invariant for the $i$-th Monge distribution. Then this invariant does not depend on the variables not involved in the subsystem.  
 \end{proposition}
  {\it Proof:} The claim follows by induction from formula (\ref{daxi}).\hfill $\Box$\\
  
 Now let us consider Darboux integrability  for the prolonged system of commuting flows. The system (\ref{commuting}) prolonged to the $r$-th order is living in $\mathbb{R}^{(2+r)n}$ with coordinates $u^i,\mu^i$ and $p^i_m:=D^m_i(\mu ^i)$,  $i=1,...,n$, $m=1,...,r $, where $D_i$ is the total derivative with respect to $u_i$. Written as an exterior differential system, (\ref{commuting}) has the form
\begin{equation}\label{extcommutingprolong}
\sigma_{i,m}:=p^i_{m+1}du^i+\sum_{k\neq i}Q^i_{m,k}du^k-dp^i_m=0,\ \ \ m=1,...,r-1,
\end{equation}
 where  one finds  $Q^i_{m,k}$ inductively by differentiation $Q^i_{m+1,k}=D_i(Q^i_{m,k})$  and by substitution of already found expressions.  We observe that

  1) the prolongation ensures that $d\sigma_{i,m}=0,$ $m=1,...,r-2, i=1,..,n$ modulo $\sigma_{i,m}=0$ and

   2) since the system (\ref{commuting}) is in involution, any its prolongation also is.\\
  Therefore it is enough to consider the exterior derivatives of $\sigma_{i,r-1}$.   Due to involution, these derivatives have the following form
$$
d\sigma_{i,r-1}=(dp^i_{r}-\sum_{k\ne i}S^i_{r,k}du^k)\wedge du^i,
$$
 where $S^i_{r,k}$ is obtained by total differentiation of  $Q^i_{r-1,k}$ with respect to $u^i$ and taking into account the 1-forms of the system.

Thus, the $i$-th Monge distribution of the prolonged system is specified by the conditions $du^i=0$ (as expected) and $dp^i_{r}=\sum_{k\ne i}S^i_{r,k}du^k$.

In what follows, we will use the notation $u=(u^1,...,u^n)$, $\mu=(\mu^1,...,\mu^n)$, $p_m=(p^1_m,...,u^n_m)$.

\begin{lemma}
If the $i$-th Monge distribution of the prolonged system  has an extra Riemann invariant of $r$-th order then this invariant can depend only on $u,\mu^i$ and $p^i_1,...,p^i_r$.
\end{lemma}
{\it Proof:}  
An extra Riemann invariant of $r$-th order, say $I(u,\mu,p_1,p_2,...p_r)$, is constant on the hyperplanes $u^i=const$ for any solution:
\begin{equation}\label{DjI}
\bar{D_j}(I)=0,\ \ \ j\ne i,
\end{equation}
where
$$
 \bar{D_j}=\partial_j+ p^j_1\partial_{\mu^j}+\sum_{s\ne j}a_{sj}(\mu^j-\mu^s)\partial_{\mu^s}+\sum_{m=1}^r\left(p^j_{m+1}\partial_{p^j_m}+\sum_{s\ne j}Q^s_{m,j}\partial_{p^s_m}\right)
$$
is the total differentiation by $u^j$ on solutions. Namely,
$$Q^s_{m,j}=D_j(p^s_m)=D_jD^m_s(\mu ^s)=D^m_sD_j(\mu ^s)=D^m_s(a_{sj}(\mu^j-\mu^s)).$$
 One easily proves by induction that  
$$D^m_s(a_{sj}(\mu^j-\mu^s))=L^m_{sj}(\mu^j-\mu^s)-\sum_{t=1}^mM^m_{tsj}p^s_t,$$
where $L^m_{sj}$ and $M^m_{tsj}$ depend only on $u$. 
Now  one successively splits equation (\ref{DjI}) with respect to $p^j_{r+1},p^j_{r},...,p^j_{1}$ and obtains $\partial_{p^j_r}(I)=\partial_{p^j_{r-1}}(I)=...=\partial_{p^j_{1}}(I)=\partial_{\mu^j}(I)=0.$
 \hfill $\Box$\\
 
 We will need a bit more detailed information on the generators of the $i$-th Monge characteristic distribution. As follows from the above proof, the vector fields $\partial_{p^j_r},\partial_{p^j_{r-1}},...,\partial_{p^j_{1}},\partial_{\mu^j},$ where $j\ne i$ belong to the distribution. Therefore one can replace the operators  $\bar{D_j}$ by 
 $$
  \tilde{D_j}=\partial_j+a_{ij}(\mu^j-\mu^i)\partial_{\mu^i}+\sum_{m=1}^r\left(L^m_{ij}(\mu^j-\mu^i)-\sum_{t=1}^mM^m_{tij}p^i_t\right)\partial_{p^i_m}.
 $$
Finally, computing the commutators 
$$
R_j:=[\partial_{\mu_j}, \tilde{D_j}]=a_{ij}\partial_{\mu^i}+\sum_{m=1}^rL^m_{ij}\partial_{p^i_m},
$$ we conclude that  the generators can be chosen as follows:
\begin{equation}\label{generatorsCom}
R_j, S_j, \partial_{p^j_r},\partial_{p^j_{r-1}},...,\partial_{p^j_{1}},\partial_{\mu^j},  \  \  \mbox{where $j\ne i$ }
\end{equation}
and 
$$
S_j=\partial_j-\sum_{m=1}^r\left(\sum_t^mM^m_{tij}p^i_t\right)\partial_{p^i_m}. 
$$
Recall that this does not mean that the span of these vector fields is the distribution but that these fields generate the distribution by inductive computation of "new" generators as Lie brackets of "old" ones  till the dimension of the span stabilizes. 
\begin{theorem} 
The $i$-th Monge distribution of  system (\ref{commuting}) has an extra Riemann invariant of minimal positive order $r$ if and only if any its $i$-th sequence of reduced Laplace transformation terminates after $r$ steps.
\end{theorem}
 {\it Proof:}  Let us apply Laplace transformation.  From  (\ref{RLaplaceMu})) we have

$$
p^i_1=D_i(\mu^i)=\left(\frac{\partial_i a_{ij}}{a_{ij}}-a_{ji}\right)\left(\bar{\mu}^i-\bar{\mu}^j\right).
$$
Therefore
$$
p^i_m=D_i^{m-1}\left[\left(\frac{\partial_i a_{ij}}{a_{ij}}-a_{ji}\right)\left(\bar{\mu}^i-\bar{\mu}^j\right)\right]
$$
is of order $m-1$ in $\bar{\mu}$. Thus, if the $i$-th Monge distribution of the prolonged system  has  an  extra Riemann invariant $I(u,\mu^i,p^i_1,p^i_2,...p^i_r)$ of order $r$ this invariant  is also the Riemann invariant of order $r-1$ for the transformed system. Now the theorem claim follows by induction from Theorem \ref{terminates1}.
 \hfill $\Box$

 \begin{theorem}
The $i$-th Monge distribution of a semihamiltonian system (\ref{hydro}) has  an extra Riemann invariant of  order $r$ if and only if the corresponding system for commuting flow (\ref{commuting}) has an extra Riemann invariant of  order $r$.
\end{theorem} {\it Proof:} 
Generic solutions to (\ref{hydro}) are related to solutions of the corresponding system (\ref{commuting}) via  (\ref{Tsarevform}). Thus a generic solution to  (\ref{hydro}) defines a so-called hydrodynamic surface in $\mathbb R^ n$. This surface can be written in an implicit form by excluding $t=\frac{\mu^j-\mu^k}{\lambda^j-\lambda^k}$ and $x=\frac{\lambda^j \mu^k-\lambda^k \mu^j}{\lambda^j-\lambda^k}$.  
One computes relations between derivatives $u^i_x$ and $\mu^i_i$ applying the total derivative by $x$ to  (\ref{Tsarevform}) (see also \cite{T-90}):
\begin{equation}\label{uxmui}
1=\sum_j(\mu^i_j-t\lambda^i_j)u^j_x=(\mu^i_i-t\lambda^i_i)u^i_x+\sum_{j\ne i} \left(a_{ij}(\mu^j-\mu^i)-\frac{\mu^j-\mu^i}{\lambda^j-\lambda^i}\lambda^i_j\right)u^j_x=\left(\mu^i_i-t\lambda^i_i\right)u^i_x.
\end{equation}

Applying the total differentiation by $x$ to this formula, one concludes inductively that $p^i_r$ can be expressed as a rational function of the following variables: 1)  $u^i_s$, $s\le r$, 2)  $u^j_q,$ where $ j\ne i,\ q<r$,  3) $t$, $\mu$, and derivatives of $\lambda$-s by $u$. Together with $\mu^s=x+\lambda^s t$, these expressions defines a map from the space of the prolonged exterior system for  (\ref{hydro}) to the space of the prolonged exterior system for  (\ref{commuting}): $(x,t,u,u_1,...,u_r)\mapsto (u,\mu,p_1,...,p_r)$. The image of this map is a submanifold $H$  obtained by excluding $x,t$ from (\ref{Tsarevform}). It is given by the following equations: 
\begin{equation}\label{H}
\mu^i(\lambda^j-\lambda^k)+\mu^j(\lambda^k-\lambda^i)+\mu^k(\lambda^i-\lambda^j)=0, \ \ i\ne j \ne k \ne i. 
\end{equation} 
Taking into account the way to obtain relations between $u_r$ and $p_r$ by differentiating the relation $$
1=(p^k_1-t\lambda^k_k)u^k_x=\left(p^k_1-\frac{\mu^k-\mu^i}{\lambda^k-\lambda^i}\lambda^k_k\right)u^i_x, 
$$
one concludes that the Monge characteristic vector of (\ref{hydro}), which is the restriction of $D_i-\lambda^iD_x$, is mapped to 
\begin{equation}\label{charH}
\zeta_i=\sum_s(\lambda^s-\lambda^i)u^s_x\bar{D}_s=\sum_s\frac{(\lambda^s-\lambda^i)^2}{p^s_1(\lambda^s-\lambda^i)-\lambda^s_s(\mu^s-\mu^i)}\bar{D}_s.
\end{equation}
Therefore any Riemann invariant of order $r$ of (\ref{commuting}), restricted on the hydrodynamic surface,  can be expressed in terms of $t,x,u,u_s$, where  $s\le r$, and of derivatives of $\lambda$-s by $u$.

Since this restriction  is constant along the lines of intersection of the hyperplanes $u^i=const$ and of the hydrodynamic surface, it is an extra Riemann invariant of order $r$ of (\ref{hydro}). It is non-trivial since it essentially depends on $u^i_r$. 

Now suppose that system (\ref{hydro}) has an extra Riemann invariant of order $r$. Since it does not depend on $u^j_r$ the image of the $i-th$ distribution of (\ref{hydro}) on $H$ is generated  by 
$$
\zeta_i,\partial_{p^j_r}, \ \ \mbox {for all $j\ne i$}.
$$
Computing $[\partial_{p^j_r},\zeta_j]$ we obtain new generators $\partial_{p^j_{r-1}}$ and inductively all $\partial_{p^j_m}$ with $j \ne i, 1\le m\le r$. Now we can replace $\zeta_i$ by 
$$
\hat{\zeta_i}=\sum_s\frac{(\lambda^s-\lambda^i)^2}{p^s_1(\lambda^s-\lambda^i)-\lambda^s_s(\mu^s-\mu^i)}\hat{D}_s,
$$
 where 
 $$
  \hat{D}_s=\partial_s+p^s_1\partial_{\mu^s}+\sum_{k\ne s} a_{ks}(\mu^s-\mu^k)\partial_{\mu^k}+\sum_{m=1}^r\left(L^m_{is}(\mu^s-\mu^i)-\sum_{t=1}^mM^m_{tis}p^i_t\right)\partial_{p^i_m}.
 $$
Further, the operators  
$$
\check{D}_j=-\frac{(p^j_1(\lambda^j-\lambda^i)-\lambda^j_j(\mu^j-\mu^i))^2}{(\lambda^j-\lambda^i)^3}[\partial_{p^j_1},\hat{\zeta_i}]
= \hat{D}_j-\frac{p^j_1(\lambda^j-\lambda^i)-\lambda^j_j(\mu^j-\mu^i)}{\lambda^j-\lambda^i}\partial_{\mu^j}
$$
are also in the distribution for $j\ne i$. 
Using them we also obtain 
$$
T_i=\hat{\zeta}_i-\sum_s\frac{(\lambda^s-\lambda^i)^2}{p^s_1(\lambda^s-\lambda^i)-\lambda^s_s(\mu^s-\mu^i)}\left(\hat{D}_s  
-\frac{p^s_1(\lambda^s-\lambda^i)-\lambda^s_s(\mu^s-\mu^i)}{\lambda^s-\lambda^i}\partial_{\mu^s}
  \right)=\sum_s(\lambda^s-\lambda^i)\partial_{\mu^s}.
$$
This operator is crucial for "splitting" $\check{D}_j$ since after computing
$$
\frac{1}{\lambda^j-\lambda^i}[T_i,\check{D}_j]=R_j+a_{ij}\sum_{s\ne i}\partial_{\mu^s}\
$$
we can replace  $\check{D}_j$ by 
\begin{equation*}
\begin{array}{l}
\bar{S_j}:=\check{D}_j-(\mu^j-\mu^i)\left(R_j+a_{ij}\sum_{s\ne i}\partial_{\mu^s}\right)=\\
\\
\ \ \ \ S_j+\frac{\lambda^j_j(\mu^j-\mu^i)}{\lambda^j-\lambda^i}\partial_{\mu^j}+\sum_{k\ne i,j} a_{kj}(\mu^j-\mu^k)\partial_{\mu^k}-a_{ij}(\mu^j-\mu^i)\sum_{s\ne i}\partial_{\mu^s}.
\end{array}
\end{equation*}
Thus, the image of the $i$-th Monge distribution of (\ref{hydro}) on $H$ is generated 
by the operators 
\begin{equation}\label{generatorsH}
\hat{\zeta_i},\ \ T_i, \ \  \bar{S}_j,\ \  R_j+a_{ij}\sum_{s\ne i}\partial_{\mu^s},\ \ \partial_{p^j_m}, \ \ j\ne i,\ \ 1\le m\le r,
\end{equation}
and has codimension at least 2 on $H$. Moreover, holds 
$$
\hat{\zeta_i}\equiv T_i\ \ \mbox{mod} \left<S_j,R_j+a_{ij}\sum_{s\ne i}\partial_{\mu^s}\right>_{j\ne i},
$$
and therefore we can exclude $\hat{\zeta_i}$ from the set of generators. Now we consider these operators as vector fields defined on the whole space of the system (\ref{commuting}) and observe that the vector fields $\partial_{\mu^j},\ \ j\ne i$ are transverse to $H$ (i.e. generate the whole tangent space with the tangent space to $H$), and that 
$$
[\partial_{\mu^k},S_j]=[\partial_{\mu^k},R_j]=[\partial_{\mu^k},T_i]=0,\  \  \ [\partial_{\mu^k},\bar{S}_j]\equiv 0,\ \ \mbox{mod} \left<\partial_{\mu^s}\right>_{s\ne i}, \forall  j,k\ne i. 
$$
Therefore the vector fields (\ref{generatorsH}) supplemented with $\left<\partial_{\mu^s}\right>_{s\ne i}$ generate a distribution of codimension at least 2 on the space of the system (\ref{commuting}). Using the vector fields $\left<\partial_{\mu^s}\right>_{s\ne i}$  we replace $\bar{S}_j,  R_j+a_{ij}\sum_{s\ne i}\partial_{\mu^s}$ be $S_j,R_j$ and obtain the set 
 (\ref{generatorsCom}). 
 \hfill $\Box$\\ 
 
 \begin{corollary}
 If  $r$-th prolongation of system (\ref{hydro}) is Darboux integrable then the  $r$-th prolongation of any commuting system is also Darboux integrable.
 \end{corollary}

\noindent {\bf Remark.} Applying Laplace transformations, one easily generalizes  Theorem \ref{commquadr} to higher orders: if a prolongation of system (\ref{commuting}) is Darboux integrable then it is integrable in quadratures. In fact, if the $i$-th Monge distribution admits an extra Riemann invariant of order $r$ then there is an $i$-sequence of Darboux transformations (possibly reduced) $\mu^i\to \mu^i_{(1)}\to ...\to \mu^i_{(q)}\to \mu^i_{(q+1)}\to...\to  \mu^i_{(r)}$: 
$$
\mu^i_{(q+1)}=-\frac{1}{b_{ij_i}^{(q)}}\partial_i\mu^i_{(q)}+\mu^i_{(q)},
$$   
and $\partial_k\mu^i_{(r)}=0.$ Thus $\mu^i_{(r)}=\varphi^i(u^i)$ and along any line parallel to $u^i$-coordinate axis one  
successively find $\mu^i_{(r-1)},\mu^i_{(r-2)},...,\mu^i$ by quadratures. The rest of the proof goes as for Theorem \ref{commquadr}.    
\section{Line congruence geometry of semihamiltonian systems with terminating sequences of Laplace transformations}\label{sectionCongr}
Semihamiltonian systems (\ref{hydro}) admit infinitely many conservation laws  \cite{T-90}, i.e.
differential one-forms $N(u)dx+M(u)dt$, closed on any solution. One can choose $n$ functionally independent densities $N^i(u)$ and rewrite (\ref{hydro}) as a system of conservation laws
(\ref{CL}).

A conservation low $N(u)dx+M(u)dt$ of (\ref{hydro}) satisfies the system
$$
\partial_iM=\lambda^i\partial_i N,~~~i=1,...,n,\ \ \ \mbox{(no summation)}.
$$
Excluding the flow  $M$, we have a consistent system of linear PDEs for the density $N$:
\begin{equation}\label{systemforN}
N_{ij}=a_{ij}N_i+a_{ji}N_j.
\end{equation}
Thus we can apply Laplace transformations to this system. These transformations are liftable to the characteristic speeds $\lambda^s$   \cite{AF-99} and give  formulas (\ref{LaplaceMu}) that we obtained earlier starting from (\ref{commuting}).

Laplace transformations for (\ref{hydro}) can be defined geometrically via the corresponding line  congruence \cite{AF-99}.
For any choice of $n$ independent conservation laws, a line of the congruence (\ref{congrn}), corresponding to a semihamiltonian system (\ref{hydro}), is tangent to $n$ {\it focal submanifolds} \cite{AF-96}, labeled by $i=1,...,n$:
$$
\begin{array}{l}
Y^0=-\lambda^{i}Y^{n+1}\\
Y^m=(M^k(u)-\lambda^{i} N^k(u))Y^{n+1},\ \ \ \ k=1,..., n.
\end{array}
$$
Passing to the affine chart $y^k=-\frac{Y^k}{Y^{n+1}}$, we get a  parametrization $r^i(u)$ of the $i$-th focal submanifold $F_i$:
\begin{equation}
\begin{array}{l}
y^0=\lambda^{i},\\
y^k=\lambda^{i} N^k(u)-M^k(u),\ \ \ \ k=1,..., n.
\end{array}
\label{focalaff}
\end{equation}
Now fix $i,j$ and consider $n$-parameter family of lines passing trough $r^i(u)$ in the direction of $\partial_j r^i(u)$:
\begin{equation}\label{paramLaplace}
\begin{array}{l}
y^0=\lambda^{i}+s\lambda^i_j,\\
y^k=\lambda^{i} N^k-M^k+s[\lambda^i_jN^k+(\lambda^i-\lambda^j)\partial_j N^k  ], \ \ \ k=1,..., n.
\end{array}
\end{equation}
Excluding the parameter $s$,  we get a line congruence corresponding to the system of conservation laws $\bar{N}^k(u)dx+\bar{M}(u)^kdt$ with
\begin{equation}\label{LaplaceCongr}
\bar{N}^k(u)=N^k-\frac{\partial_jN^k}{a_{ij}},\ \ \ \  \bar{M}^k=M^k-\frac{\lambda^i\partial_j N^k}{a_{ij}}.
\end{equation}
It is quite natural to call (\ref{LaplaceCongr}) {\it Laplace transformation} of the congruence (\ref{congrn}).

For any $i\ne k\ne j\ne i$ we check by direct computation, using (\ref{systemforN}) and (\ref{compatibilityCom}), that
$$
\partial_k \bar{M}=\frac{a_{ij}\lambda^k-a_{ki}\lambda^i}{a_{ij}-a_{ki}}\partial_k\bar{N}.
$$
Similarly, (\ref{systemforN}) implies
$$
\partial_i \bar{M}=\left(\lambda^i-\frac{\lambda^i_i}{ \frac{\partial_i a_{ij}}{a_{ij}}-a_{ji}    }\right)\partial_i\bar{N}
$$
and
$$
\partial_j \bar{M}=\lambda^i\partial_j\bar{N}.
$$
Thus the characteristic speeds of the transformed system of conservation laws {\it do not depend on the choice of $N^i$} and coincide with ones of (\ref{hydro}) after applying the Laplace transformation (\ref{LaplaceMu}) to $\lambda$-s.

 \medskip
\noindent {\bf Remark 1.}
In fact, it is not difficult to see that semihamiltonian systems can be defined as strictly hyperbolic systems of conservation laws admitting $n$ Riemann invariants \cite{S-94}: the compatibility conditions (\ref{compatibilityCom}) for the commuting flows are also the compatibility conditions for the system  (\ref{systemforN}) for the densities of conservation laws for the diagonal system (\ref{hydro}).

\medskip

Let us translate Darboux integrability into the geometric language.
To this end we recall the notion of {\it reciprocal transformation} for the system (\ref{hydro}): given 2 conservation laws $Bdx+Adt$, $Ndx+Mdt$, we introduce new independent variables $X,T$ by
\begin{equation}
\begin{array}{c}
dX=B(u)dx+A(u)dt,\\
dT=N(u)dx+M(u)dt
\end{array}
\label{Trec}
\end{equation}
and rewrite (\ref{hydro}) as
\begin{equation}\label{hydroT}
u^i_T=\left(\frac{B\lambda^i - A}{M-N\lambda^i}\right)u^i_X.
\end{equation}
Reciprocal transformations respect the property of being semihamiltonian, but not the property of being Darboux integrable. One easily sees this already for extra Riemann invariants of order 0: the condition $\lambda^i_k=0,\ \forall k\ne i$ does not imply $\partial_k \left(\frac{B\lambda^i - A}{M-N\lambda^i}\right)=0, \ \forall k\ne i$. The reason is that  an extra Riemann invariant depends on $x,t$, which became nonlocal variables (i.e. recoverable via integration of forms closed on solutions) after applying a nontrivial reciprocal transformation. If one uses only linear combinations of $dx,dt$ and of  the "coordinate" conservation laws  $N^i(u)dx+M^i(u)dt$ then the reciprocal transformations form the projective group $PGL_{n+2}$ acting in the projective space of the congruence (\ref{congrn}) \cite{AF-96}.

To understand the geometry of Darboux integrability it is enough to interpret, for a fixed $i$, the conditions $a_{im}=0,$ $\forall \ m\ne i$ or $b_{im}=0,$ $\forall \ m\ne i$  since it is the form of the system after a (sequence of) Laplace transformations of the congruence.

A line of the congruence (\ref{congrn}) is tangent to the $i$-th focal submanifold $F_i$, parametrized by (\ref{focalaff}), at some point $r^i(u)=(\lambda^{i},
\lambda^{i} N^1-M^1, ..., \lambda^{i} N^n(u)-M^n)$. If $u$ moves along the $k$-th coordinate line, given by $u^l=const,\ \forall l\ne k$, the point $r^i(u)$ traces  a curve on the focal submanifold. Repeating these construction for all $r^i(u)$ and all $k$, we obtain the so-called
{\it focal net}. Note that the congruence lines are tangent to the images of the $i$-th coordinate lines on the $i$-th focal submanifold $F_i$. The $ij$-Laplace transformation replaces the congruence lines by the lines tangent to the images of the $j$-th coordinate lines. The crucial fact about diagonalizable systems of conservation laws is that the curves of focal net are pairwise conjugate \cite{AF-99}, which translates analytically as
$$
\partial_k\partial_l r^i(u)\in T_{r^i}F_i.
$$
Thus the families of tangent planes to $F_i$ along the images of  $i$-th coordinate lines envelope a ruled hypersurfaces, whose $(n-1)$-dimensional generators are spanned  by the vectors $\partial_k r^i(u),\  k\ne i$. These generators determine an integrable $(n-1)$-dimensional distribution on $F_i$. Let as call the integral submanifolds of this distribution {\it conjugate to the congruence with respect to the $i$-th focal submanifold}.
\begin{theorem}\label{aimRec}
Suppose that, for a  semihamiltonian system (\ref{hydro}),  the $i$-th focal submanifold of the corresponding congruence is $n$-dimensional.   If the $i$-th Monge characteristic distribution  has an extra Riemann invariant of order 0 then the  submaifolds, conjugate to the congruence with respect to the $i$-th focal submanifold $F_i$, are cut on $F_i$ by a pencil of hyperplanes.

Conversely, if for  some conservative form (\ref{CL}) of (\ref{hydro}) the  submaifolds, conjugate to the congruence with respect to the $i$-th focal submanifold $F_i$, are cut on $F_i$ by a pencil of hyperplanes then  there is a reciprocal transformation of (\ref{hydro}), whose $i$-th Monge characteristic distribution  has an extra Riemann invariant of order 0.
\end{theorem}
 {\it Proof:}   If there is a Riemann invariant of order 0, then $\partial_k\lambda ^i=0, \forall k\ne i$ and the conjugate submanifolds are cut on $F_i$, parametrized as in (\ref{focalaff}), by the pencil of hyperplanes $y^0=\lambda^i(u).$

To prove the converse claim, we choose homogeneous coordinates so that the pencil of hyperplanes are spanned by   $Y^0=0$ and $Y^{n+1}=0.$  Changing the homogeneous coordinates corresponds to applying a reciprocal transformation to (\ref{hydro}).
 \hfill $\Box$

\medskip

\noindent {\bf Remark 2.} The conjugacy of the directions along $k$-th and $l$-th coordinate lines can be written as $$\partial_k\partial_l r^i(u)=\sum_{m=1}^nC_{im}(u)\partial_m r^i(u).$$ For $n=2$ the congruence sits in $\mathbb{P}^3$ and all components $y(u)$ of the parametrization   (\ref{focalaff}) of the focal surface $F_i$ satisfy a scalar equation of the form (\ref{2ndOrder}). There is a nice geometric interpretation of terminating  in one step of Laplace  transformation {\it of the corresponding system (\ref{2ndOrder})}: the  focal surfaces of the congruence, obtained by Laplace transformation, degenerate to curves (see \cite{D-96}, V2, book 4, cap.I). This interpretation is not valid for the Laplace transformation of (\ref{hydro}). The reason is that the equation for $\lambda $ is {\it different} from the one for {\it all} the coordinates of the corresponding focal surface, therefore  the Laplace transformation, given geometrically by the same construction, is applied to different systems!

\medskip

 Vanishing of $b_{ij} \  \forall j\ne i$ for some fixed $i$ implies existence of Riemann invariant of order 1 (see lemma \ref{order01}) for the $i$-th Monge distribution. This  case has the following geometric interpretation.

\begin{theorem}
The envelopes of tangent planes to the $i$-th focal hypersurface along $i$-th coordinate lines $u^k=const\ \forall k\ne i$ are cones over the images of these lines on the focal hypersurface with a $(n-3)$-dimensional vertex if and only if $b_{ij}=0 \ \forall j\ne i$.
\end{theorem}
 {\it Proof:} The envelopes are swept by $(n-1)$-dimensional affine planes, generated by vectors  $\vec{v}_j-a_{ij}\vec{v}_0, \ {k\ne i}$ at $r^i(u)=(\lambda^{i},
\lambda^{i} N^1-M^1, ..., \lambda^{i} N^n(u)-M^n)$, where $\vec{v}_0=(1,N^1,...,N^n)$, $\vec{v}_j=\partial_j\vec{v}_0$. Using (\ref{systemforN}), one checks that $\partial_i(\vec{v}_j-a_{ij}\vec{v}_0)\in \langle \vec{v}_k-a_{ik}\vec{v}_0 \rangle_{k\ne i}$ if and only if $\partial_i a_{ij}-a_{ji}a_{ij}=a_{ij}b_{ij}=0$. This means that the vector space $\langle \vec{v}_k-a_{ik}\vec{v}_0 \rangle_{k\ne i}$ is stable along the $i$-th coordinate line $u^k=const,\ k\ne i$ and all the affine generators intersect at some $(n-3)$-dimensional space at the infinite hyperplane.
\hfill $\Box$

\section{Subalgebras of commuting flows and Darboux integrability}

Let as fix $i$ and look for subalgebras of the commuting flow algebra (\ref{commuting}) cut by a linear differential relation of the form
\begin{equation}\label{diffrelation}
\mu^i_i=\sum_{k\ne i}K_k(u)\mu^k_k+\sum_jL_j(u)\mu^j,
\end{equation}
compatible with (\ref{commuting}).  Note that  up to index permutation this is a generic form of one linear (homogeneous) differential relation of first order.

\begin{theorem}\label{subalgebra1}
Equations (\ref{commuting}) for commuting flows of a semihamiltonian system (\ref{hydro}) are compatible with the differential relation (\ref{diffrelation}) if and only if $K_k=0\ \forall k$, $L_j=0,\ \forall  j\ne i$, $L_i=b_{ij}\ \forall  j\ne i$ and the $i$-th Monge characteristic distribution admits an extra Riemann invariant of order 1, i.e. $\partial_jb_{ij}+a_{ij}b_{ij}=0 \ \ \forall \ j\ne i$.
\end{theorem}
 {\it Proof:}
For any $j\ne i$ we can compute the second derivatives $\mu^i_{ij}$ in two different ways: either differentiating (\ref{diffrelation}) by $u^j$ or differentiating $a_{ij}(\mu^j-\mu^i)$ by $u^i$. Simplifying modulo (\ref{commuting},\ref{diffrelation}), we see that $\mu^i_{ii}$  appears only once with the coefficient $K_j$. Therefore $K_j=0$, since we are not allowed to get any further relations for (\ref{diffrelation}) compatible with (\ref{commuting}). Similarly, the coefficient of $\mu^j_j$ yields  $L_j=0$, vanishing of the  coefficient of $\mu^j$ implies $L_i=b_{ij}$ , and, finally,  on  computing the coefficient of $\mu^i$ we conclude $\partial_j L_{i}=-a_{ij}+a_{ij}a_{ji}$, hence  $\partial_jb_{ij}+a_{ij}b_{ij}=0$.
 \hfill $\Box$\\

Even though the infinitesimal action of the symmetry algebra  of commuting flows is defined only on solutions of (\ref{hydro}) via (\ref{hydroMu}) (i.e. solutions to (\ref{commuting}) provide not  classical but generalized symmetry operators \cite{O-86}) we still can talk about the  orbits of a solution.
\begin{theorem}\label{orbitth}
Suppose that  the $i$-th Monge characteristic distribution  of a semihamiltonian system (\ref{hydro}) admits an extra Riemann invariant of order 1. Then the orbit of the subalgebra of (\ref{commuting}) defined by the differential relation
\begin{equation}
\mu^i_i=b_{ij}\mu^i,
\end{equation}
is given by the differential relation (\ref{orbit}) and is fixed by the choice of  $\varphi_i$.
\end{theorem}
 {\it Proof:}
  Following the general approach (see for instance \cite{O-82}), we compute a differential equation that defines the orbit of a generic solution.   Let $\mu$ be a solution to (\ref{commuting}). Then the operator
 \begin{equation}\label{genSym}
 \sum_j \mu^j u^j_x\partial _j+  \sum_j D_x(\mu^j u^j_x)\partial _{u^j_x}+\sum_j D_t(\mu^j u^j_x)\partial _{u^j_t}-\mu^iD_x
 \end{equation}
gives a generalized symmetry of (\ref{hydro}). (We used the fact that the total differentiations $D_x,D_t$ generate the ideal of generalized symmetries of {\it any} system  \cite{O-86}.)  Let us consider the action of (\ref{genSym}) on differential functions of $x,t,u,u_x$.   Substituting $\mu^j_k$ for $k\ne j$ from  (\ref{commuting})  we get
$$
-\mu^i(\partial_x+u^i_x\partial_i+u^i_{xx}\partial_{u^i_x})+\sum_{j}\mu^j(u^j_x\partial_j+u^j_{xx}\partial_{u^j_x})+\sum_ju^j_x\left(\mu^j_ju^j_x+\sum_{k\ne j}a_{jk}(\mu^k-\mu^j)u^k_x \right)\partial_{u^j_x}.
$$
 Collecting the (operator-valued) coefficients of $\mu^j$, $\mu^j_j$ and substituting $\mu^i_i$ yields:
 \begin{eqnarray}
 \mu^i\left(-\partial_x+b_{ij}(u^i_x)^2\partial_{u^i_x}-\sum_{j\ne i}u^j_x\left(\partial_j+a_{ij}u^i_x\partial_{u^i_x}\right)\right)+\\\sum_{j\ne i}\mu^j_j(u^j_x)^2\partial_{u^j_x}+
 \sum_{j\ne i}\mu^j\left(u^j_x\left(\partial_j+a_{ij}u^i_x\partial_{u^i_x}\right)  +u^j_{xx}\partial_{u^j_x}\right).
 \end{eqnarray}
To write down an equation for the orbit, we need invariants of  our  subalgebra, depending on $x,t,u,u_x$. On the "tangent space" of the orbit at a "point" (which is a solution to (\ref{hydro})), the initial values of $\mu^k$, $k=1,...,n$ and $\mu^k_k$, $k\ne i$ are arbitrary. Thus, splitting with respect to these variables, we see that the invariants are the first integrals of
the distribution generated by the following set of vector fields
\begin{equation}\label{pseudo}
\{\partial_{u^j_x}, \ \partial_j+a_{ij}u^i_x\partial_{u^i_x}, \ \ -\partial_x+b_{ik}(u^i_x)^2\partial_{u^i_x}\}_{j\ne i}.
\end{equation}
 Note that $b_{ik}=b_{ij}$ for any $j\ne i \ne k\ne j$ and that  this distribution is involutive  by Lemma \ref{order01}. A solution  to (\ref{hydro}) gives a 2-dimensional surface in $(2n+2)$-dimensional space with the coordinates $x,t,u,u_x$. If this solution is generic then its orbit, being also the orbit of the Lie pseudogroup with the $(2n-1)$-dimensional algebra (\ref{pseudo}), is $(2n+1)$ dimensional. Being invariant  of codimension 1 in  $(2n+2)$-dimensional space, the orbit is given by a single relation between invariants. Thus, one can define the orbit by a single equation
\begin{equation}\label{iorbit}
F(x,t,u,u^i_x)=0,
\end{equation}
where $F$ is an invariant of (\ref{pseudo}) and this equation is compatible with (\ref{hydro}). (For more detail on the construction of orbits, i.e. automorphic PDEs, see \cite{O-82}).

Now we use the compatibility. To this end we substitute the  differentials
$$
du^i_x=-\frac{1}{\partial_{u^i_x}F}\left(  \partial_t(F)dt + \partial_x(F)dx+\sum_k \partial_i(F)du^k  \right),\ \ \ \
du^k=u^k_x(dx+\lambda^k dt).\\
$$
into
$$
d(du^i)=du^i_x\wedge (dx+\lambda^i dt)+u^i_x\sum_k \lambda^i_ku^k_xdx\wedge dt=0.
$$
The resulting equation
$$
\{-\partial_t  +\lambda^i\partial_x-\lambda^i_i(u^i_x)^2\partial_{u^i_x}\}F+\sum_{k\ne i}u^k_x(\lambda^i-\lambda^k)\{\partial_k+a_{ik}u^i_x\partial_{u^i_x}\}F=0
$$
 cannot impose any additional relation on the orbit. Since  the operators  $\partial_k+a_{ik}u^i_x\partial_{u^i_x}$ kill $F$ for any $k\ne i$ the function $F$ is also invariant along the vector field $-\partial_t  +\lambda^i\partial_x-\lambda^i_i(u^i_x)^2\partial_{u^i_x} $. Invoking Lemma \ref{order01} we conclude that $F$ is a Riemann invariant of the $i$-th Monge distribution. Since the integral manifolds of this distribution are of codimension 2, any single relation (\ref{iorbit}) may be rewritten as (\ref{orbit}).
   \hfill $\Box$\\
   
Theorem \ref{subalgebra1} is generalizable for prolonged systems as follows. A general linear homogeneous differential  relation of order $r$ cutting a subalgebra of (\ref{commuting}) has the form
\begin{equation}\label{diffrelationr}
p^i_r=\sum_{k\ne i}K_k(u)p^k_r+\sum_{j,1\le m<r}L_{jm}(u)p^j_m+\sum_jL_j(u)\mu^j.
\end{equation}
\begin{theorem}\label{subalgebrak}
Equations (\ref{commuting}) for commuting flows of a semihamiltonian system (\ref{hydro}) are compatible with the differential relation (\ref{diffrelationr}) of minimal possible order $r$ if and only if the $i$-th Monge characteristic distribution admits an extra Riemann invariant of minimal order $r$ and the relation has the form
$$
D_i \mu^i_{(r-1)}=b^{(r-1)}_{ij}\mu^i_{(r-1)},
$$ 
where $\mu^i_{(r-1)}$ is the result of applying an $i$-sequence of Laplace transform of length $r-1$, written in terms of $\mu^i,p^i_1,...,p^i_{r-1}$, and $b^{(r-1)}_{ij}$ corresponds to the resulting system, i.e.
$K_k=0\ \forall k$, $L_j=L_{jm}=0\ \forall  j\ne i$.
\end{theorem}
 {\it Proof:} Applying the total differentiation $D_k$, $k\ne i$ to the both sides of (\ref{diffrelationr}) we observe that $D_k(p^i_r)=D_i^r(a_{ik}(\mu^k-\mu^i))$ do not contain derivatives of order $r+1$, while the r.h.s. has the term $K_kp^k_{r+1}$. Therefore $K_k=0$ for $k\ne i$. Now the claim follows by induction form Theorem  \ref{subalgebra1} after applying Laplace transform. \hfill $\Box$\\
 
In an obvious way one generalizes Theorem \ref{orbitth} to higher order  and proves the generalized Theorem by induction.  

\section{Concluding remarks}

\subsection{Reciprocal transformations}

Theorem (\ref{aimRec}) does not give a geometric criterion for a semihamiltonian systems (\ref{hydro}) to be Darboux integrable up to a reciprocal transformation: such transformations correspond to projective transformations of the corresponding congruence only for linear combinations of $dx,dt$ and of  the "basic" conservation laws  $N^i(u)dx+M^i(u)dt$. Semihamiltonian systems (\ref{hydro}) have infinitely many linearly independent conservation laws and therefore infinitely many reciprocal transformations not reducible to projective transformations of the congruence.  Reciprocal transformation can bring a system (\ref{hydro}) that is not Darboux integrable to one that possess extra Riemann invariants for some or all its Monge distributions. A criterion for Darboiux integrability by first order Riemann invariants after reciprocal transformation was found in \cite{AFN-17} though in relation to a quite different problem. This interesting question deserves a further investigation.  

\subsection{Nondiagonalizable systems} The following  hamiltonian system
$$
\frac{\partial}{\partial t}\left(\begin{array}{l} u_1\\ u_2\\ u_3 \end{array} \right)=\left(\begin{array}{ccc} -1 & 1 & 1\\ 1 & -1 & 1\\ 1 & 1 & -1 \end{array} \right)\frac{d}{dx}\left(\begin{array}{l} \partial_{u_1} H(u)\\  \partial_{u_2} H(u)\\  \partial_{u_3} H(u) \end{array} \right)
$$
with the Hamiltonian density $H(u)=\frac{1}{2}u_1u_2u_3$ is nondiagonalizable, linearly degenerate and "integrable" in the sense of having infinitely many higher conservation laws and admitting a Lax pair. It does not possesses Riemann invariants of order 0 or 1. Its characteristic speeds can be made constant by a  reciprocal transformation \cite{F-93}, thus giving Riemann invariants of order 0 for the transformed system. The author does not know if this system has Riemann invariants of higher order. It seems that no example of Darboux integrable non-diagonalizable system of hydrodynamic type has ever appeared in literature.

\section*{Acknowledgements}
This research was supported by FAPESP grants \#2018/20009-6 and \#2021/10380-1.

\end{document}